\documentclass{article}
\usepackage{amsmath}
\usepackage{amssymb}
\usepackage{amsfonts}

\newcommand{\sn}{\smallskip\noindent}
\newcommand{\mn}{\medskip\noindent}
\newcommand{\bn}{\bigskip\noindent}

\newcommand{\qed}{\vrule height 1.2ex width 1.1ex depth -.1ex}

\newcommand{\A}{\mathcal{A}}
\newcommand{\B}{\mathcal{B}}
\newcommand{\R}{\mathcal{R}}
\newcommand{\Hh}{\mathcal{H}}
\newcommand{\Oo}{\mathcal{O}}
\newcommand{\G}{\mathcal{G}}
\newcommand{\D}{\mathcal{D}}
\newcommand{\Z}{\mathcal{Z}}
\newcommand{\E}{\mathcal{E}}
\newcommand{\T}{\mathcal{T}}

\newcommand{\Ss}{\mathcal{S}}
\newcommand{\dd}{\mathrm{d}}
\newcommand{\id}{\mathrm{id}}
\newcommand{\kere}{\mathrm{ker}}
\newcommand{\Lin}{\mathrm{Lin}}
\newcommand{\bbbn}{\mathbb{N}}
\newcommand{\bbbz}{\mathbb{Z}}
\newcommand{\bbbr}{\mathbb{R}}

\newcommand{\dif}{\mathrm{d}}

\author{Konrad Schm\"udgen}
\title{\bf Commutator Representations of Differential Calculi on the 
Quantum Group ${\bf SU_q(2)}$}
\date{\small \sc Universit\"at Leipzig,
Fakult\"at f\"ur Mathematik und Informatik,
Augustusplatz 10/11, D-04109 Leipzig, Germany\\
E-mail: schmuedg@mathematik.uni-leipzig.de}

\begin{document}

\maketitle

\begin{abstract}\noindent
Let $(\Gamma,\dd)$ be the $3D$-calculus or the $4D_\pm$-calculus on the 
quantum group $SU_q(2)$. We describe all pairs $(\pi,F)$ of a 
$\ast$-representation $\pi$ of $\Oo(SU_q(2))$ and of a symmetric 
operator $F$ on the representation space satisfying a technical condition 
concerning its domain such that there exist a homomorphism of first order 
differential calculi which maps $\dd x$ into the commutator 
$[{\rm i}F,\pi (x)]$ for $x\in\Oo(SU_q(2))$. As an application 
commmutator representations of the 2-dimensional left-covariant
calculus on Podles quantum 2-sphere $S^2_{qc}$ with $c=0$ are given.
\end{abstract}
Subj. Class.: Quantum groups, Non-commutative geometry
\mn

\noindent
{\bf 1. Introduction}

\mn
There are various ways to develop a noncommutative differential 
calculus on a given $\ast$-algebra $\A$. Some approaches are based 
on derivations [3,7] of the algebra $\A$, while others 
consider differential forms as the basic objects. In Alain 
Connes' programm [2] of noncommutative geometry, the fundamental concept 
for the quantized calculus is the $K$-cycle. Other 
names for closely related notions are Fredholm modules and spectral triples. 
If we omit technical subtleties, then the underlying idea of all these 
concepts is easy to explain: One has a $\ast$-representation $\pi$ of the 
algebra $\A$ and a self-adjoint operator $F$ on the representation space, 
and the differentiation of an algebra element $a$ is defined by the 
commutator of the operators i$F$ and $\pi(a)$, where i is the imaginary unit. 

On the other hand, for quantum groups there is a well-developed theory of 
covariant differential calculi which was initiated by S.L. Woronowicz [15], 
see also the 
monograph [6]. The relations between this theory and Connes' approach 
via $K$-cycles or spectral triples are still open. Let us be more 
specific and consider the simplest non-trivial compact quantum group 
$SU_q(2)$. Then we have three distinguished left-covariant differential 
calculi on the Hopf $\ast$-algebra $\Oo(SU_q(2))$: 
the $3D$-calculus and the $4D_\pm$-calculi.
 All three calculi are due to 
S.L. Woronowicz [14],[15]. The $4D_\pm$-calculi are 
even bicovariant, but we shall not use this here. 
These calculi have a number of nice 
properties, and there is a common belief that they are favoured 
candidates for the study of noncommutative geometry on the quantum 
group $SU_q(2)$. Thus it seems to be natural to ask whether or not 
one of these calculi can be described by means of a $K$-cycle or 
more generally by the commutators $[\mbox{i} F,\pi (\cdot)]$ for some 
$\ast$-representation $\pi$ of $\Oo(SU_q(2))$ and some 
symmetric operator $F$.

Under an additional technical condition called admissibility we 
describe all such pairs $(\pi,F)$ which represent the $3D$-calculus or the 
$4D_\pm$-calculus on $SU_q(2)$. The main results about this matter 
are stated as 
Theorems 1--3 in Section 4, while the proofs of these results are postponed 
to Section 7. It turns out that the $3D$-calculus can be 
faithfully represented 
as a commutator. Some properties of such commutator representation of the
$3D$-calculus and some examples are treated in Section 5. An application 
to the Podles quantum sphere $S^2_{qc}$ with $c=0$ is sketched in Section 6. 
The $4D_{\pm}$-calculus does not admit a faithful admissible 
commutator representation,
because any such representation passes to 
a 3-dimensional left-covariant quotient calculus. Further, 
we prove that neither for the 3D-calculus nor for 
the $4D_{\pm}$-calculus there exists a non-trivial commutator 
representation $(\pi,F)$ such that all operators 
$[{\rm i}F,\pi(x)]$, $x \in \A,$ are bounded. This shows in particular 
that none of these calculi can be given by means of a K-cycle in the sense 
of A. Connes.

In Section 2 we collect some basic definitions on differential 
calculi 
and some simple facts needed later. In Section 3 
we repeat the structure of a general $\ast$-representation of 
the $\ast$-algebra $\Oo(SU_q(2))$ and some facts about the $3D$- and 
the $4D_{\pm}$-calculi on $SU_q(2)$. Further, we develop the 
$2$-dimensional calculus on the Podles sphere $S_{qc}^2$ for $c=0$ 
as the induced FODC of the $3D$-calculus on $SU_q(2)$.

Let us fix some general notation. We use the 
Sweedler notation $\Delta (x)=x_{(1)}\otimes x_{(2)}$ for the 
comultiplication of a Hopf algebra element $x$. Let 
$x(n)y$ denote the equation which is obtained by multiplying equation $(n)$ 
by $x$ from the left and by $y$ from the right. Throughout 
we set 
\begin{equation*}
\lambda := q-q^{-1} ~{\rm and}~ \lambda_+ := q+q^{-1}.
\end{equation*}

\bn
{\bf 2. Commutator representations of first order differential calculi } 

\sn
Let $\A$ be a complex $\ast$-algebra with unit element. The involution of 
$\A$ is denoted by $\ast$. 

A {\it first order differential calculus} (abbreviated, a FODC) over  $\A$ is a pair $(\Gamma,\dd)$ of an $\A$-bimodule 
$\Gamma$ with a linear mapping $\dd:\A\rightarrow \Gamma$ such 
that the following two conditions hold:
\begin{enumerate}
\item[(i)] $\dd$ satisfies the Leibniz rule $\dd(xy) = 
x{\cdot} \dd y+\dd x{\cdot}  y~{\rm for}~x,y\in\A$,
\item[(ii)] $\Gamma={\rm Lin}\{x{\cdot} \dd y{\cdot} z; x,y,z \in\A\}.$
\end{enumerate}
A {\it first order differential $\ast$-calculus} (briefly, a 
$\ast$-FODC) over  $\A$ is a FODC $(\Gamma,\dd)$ equipped with 
an involution $\ast:\Gamma\rightarrow\Gamma$ of the complex vector 
space $\Gamma$ such that
\begin{enumerate}
\item[(iii)] $(x{\cdot} \dd y{\cdot} z)^\ast=z^\ast{\cdot} 
d(y^\ast){\cdot} x^\ast ~{\rm for}~ x,y,z\in\A.$
\end{enumerate}

By a {\it homomorphism} of a FODC $(\Gamma_1,\dd_1)$ 
into a FODC $(\Gamma_2,\dd_2)$ over $\A$ we mean a 
linear mapping $\rho:\Gamma_1\rightarrow \Gamma_2$ such that 
$\rho (x{\cdot} \dd_1 y{\cdot} z)=x{\cdot} \dd_2y{\cdot} z$ 
for $x,y,z\in \A$. A {\it homomorphism} of a $\ast$-FODC 
$(\Gamma_1,\dd_1,\ast_1)$ into a $\ast$-FODC $(\Gamma_2,\dd_2, \ast_2)$ 
is a FODC homomorphism $\rho:\Gamma_1\rightarrow \Gamma_2$ such that 
$\rho (\omega^{\ast_1})=\rho (w)^{\ast_2}$ for  
$\omega\in\Gamma_1$. If no ambiguity can arise, we denote a 
FODC $(\Gamma,\dd)$ or a $\ast$-FODC 
$(\Gamma,\dd,\ast)$ by $\Gamma$.

Now let $\A$ be a Hopf $\ast$-algebra. 
A $\ast$-FODC $\Gamma$ over $\A$ is called 
{\it left-covariant} if there exists a linear mapping 
$\varphi :\Gamma\rightarrow \A\otimes \Gamma$ such that $
\varphi (x{\cdot} \dd y)=\Delta (x) (\id \otimes \dd)\Delta (y)$ 
for all $x,y\in\A$. Suppose $\Gamma$ is a left-covariant $\ast$-FODC. 
We define
$$
\omega_\Gamma (x) := S(x_{(1)})\dd x_{(2)},~ x\in\A,~{\rm and}~
R_\Gamma := \{ x\in\ker\varepsilon  :~\omega_\Gamma (x)=0\}.
$$
It is well-known [15,6] that $\R_\Gamma$ is a right ideal of $\A$ 
which characterizes the left-covariant FODC $\Gamma$ up to 
isomorphism. 

Suppose that $\pi$ is a $\ast$-representation of the $\ast$-algebra $\A$ by 
bounded operators on a Hilbert space $\Hh$ and $F$ is a symmetric linear 
operator on  $\Hh$ with dense domain $\D (F)$. Let us assume that there 
exists a linear subspace $\D_F$ of $\D(F)$ such that $\D_F$ is dense in 
$\Hh$ and $\pi(\A)\D_F\subseteq\D (F)$. Let $\Gamma_{\pi,F}$ denote the 
linear span of operators
$$ 
\pi(x) (F\pi (y)-\pi (y) F)\pi(z)\lceil\D_F,~ x,y,z\in\A,
$$
where the symbol $\lceil \D_F$ denotes the restriction of the 
corresponding operator to $\D_F$. It is clear that $\Gamma_{\pi,F}$ 
is an $\A$-bimodule with left and right action of an element $x\in\A$ 
given by multiplication by $\pi(x)$ from the left and the right, 
respectively. Define a linear mapping 
$\dd_{\pi,F}:\A\rightarrow\Gamma_{\pi,F}$ and an antilinear 
mapping $\ast:\Gamma_{\pi,F}\rightarrow\Gamma_{\pi,F}$ by
$$
\dd_{\pi,F} (x):=({\rm i}F\pi (x)-\pi (x){\rm i}F)\lceil \D_F, \\ x\in\A,
$$
where ${\rm i}$ is the imaginary unit, and
$$
T^\ast := {\sum\nolimits_j} \pi (z^\ast_j)(\pi 
(y^\ast_j)F- F\pi(y^\ast_j))\pi(x^\ast_j)\lceil \D_F,
$$
where
$$
T=\sum\nolimits_j \pi (x_j)(F\pi(y_j)-\pi(y_j)F)\pi(z_j)\lceil \D_F \in 
\Gamma_{\pi,F}.
$$
Using the facts that $\D_F$ is dense in $\Hh$, the operator $F$ is 
symmetric and $\pi$ is a $\ast$-representation it is easy to check 
that the mapping $T \to T^\ast$ is well-defined (that is, 
$T^\ast=0$ when $T=0$). Further, it is not difficult to verify that 
the triple $(\Gamma_{\pi, F}, \dd_{\pi, F},\ast)$ is a $\ast$-FODC 
over $\A$. We call it the $\ast${\it-FODC associated with the pair} 
$(\pi,F)$ and denote it simply by $\Gamma_{\pi,F}$. (With a few 
modifications concerning the domains the preceding 
construction carries over to unbounded $\ast$-representations 
$\pi$ as well. We shall not need this in this paper, because we 
consider the coordinate $\ast$-algebra $\Oo(SU_q(2))$ 
which has only bounded $\ast$-representations.) 
The above notion is closely related to Alain Connes' concept 
of a $K$-cycle$^1$. If in addition $F$ is a self-adjoint 
operator with compact resolvent and the commutator $[\pi(x),F]$ is 
bounded for any $x\in \A$ the pair $(\pi,F)$ is called a $K$-{\it cycle} 
over the $\ast$-algebra $\A$.

Now let $\Gamma$ be an arbitrary $\ast$-FODC over $\A$ and let $\pi$ and 
$F$ be  as above. We shall say that the pair $(\pi,F)$ is a 
{\it commutator representation} of the $\ast$-FODC $\Gamma$ if there 
exists a homomorphism $\rho$ of the $\ast$-FODC $\Gamma$ to the 
$\ast$-FODC $\Gamma_{\pi,F}$ associated with $(\pi, F)$. If $\rho$ is 
injective, then $(\pi,F)$ is called a {\it faithful} commutator 
representation of $\Gamma$. A slight reformulation of this definition is 
given by 

\mn
{\bf Lemma 1.} {\it A pair $(\pi,F)$ as above is a commutator representation 
of the $\ast$-FODC $\Gamma$ if and only if $\sum_jx_j\dd y_j =0$ 
in $\Gamma$ with $x_j, y_j\in\A$ always implies that 
$\sum_j\pi (x_j) ({\rm i}F\pi (y_j)-\pi (y_j){\rm i}F)\lceil \D_F=0$.}

\bn
{\bf Proof.} The only if part is trivial. If this condition is fulfilled, 
then there exists a well-defined linear map 
$\rho:\Gamma\rightarrow\Gamma_{\pi,\Gamma}$ such that 
$\rho(\sum_j x_j\dd{y_j})=\sum_j\pi(x_j)({\rm i}F\pi(y_j)-\pi(y_j){\rm i}F)
\lceil\D_F$. One easily checks that $\rho$ is a homomorphism of 
$\Gamma$ to $\Gamma_{\pi,F}$.\hfill\qed

\sn
Let $(\pi,F)$ be a pair of a $\ast$-representation $\pi$ of $\A$ and a 
symmetric linear operator as above. For $x\in\A$ we define the linear 
operator 
\begin{equation}\label{l2-1}
\Omega_{\pi,F}(x):=(\pi(S(x_{(1)})) F\pi(x_{(2)})-\varepsilon (x)F) 
\lceil \D_F.
\end{equation}

The following very simple observations are needed in the proofs of 
the main theorems below.

\bn
{\bf Lemma 2.} {\it Suppose that $(\pi,F)$ is a commutator 
representation of the left-covariant FODC $\Gamma$ and let $\rho$ 
denote a FODC homomorphism of $\Gamma$ to $\Gamma_{\pi,F}$. Then we 
have $\rho(\omega_\Gamma(x)) ={\rm i} \Omega_{\pi,F}(x)$ for all 
$x \in \A$. In particular, if $x$ belongs to the right ideal $R_\Gamma$ 
associated with $\Gamma$, then $\Omega_{\pi,F}(x)=0$.}

\mn
{\bf Proof.}  
Using the definitions of $\omega_\Gamma(x)$ and $\Omega_{\pi,F}(x)$ and 
the fact 
that $\rho$ is a FODC homomorphism we compute 
\begin{gather*}
\quad\quad \rho(\omega_\Gamma (x))=\rho (S(x_{(1)}))\rho(\dd x_{(2)})=\rho 
(S(x_{(1)}))\dd_{\pi,F}(x_{(2)})\\
\qquad\qquad\qquad\qquad\qquad=\pi(S(x_{(1)}))({\rm i}F\pi (x_{(2)})-
\pi(x_{(2)}){\rm i}F) \lceil \D_F 
={\rm i} \Omega_{\pi,F} (x).\qquad\qquad\qquad\qquad
\hfill\qed
\end{gather*}

\mn
{\bf Lemma 3.} {\it Let $\rho: \Gamma_1 \to \Gamma_2$ be a 
homomorphism of the left-covariant FODC $\Gamma_1$ into another 
FODC $\Gamma_2$. Suppose 
that $\B$ is a subset of ${\rm ker}~ \varepsilon$ such that $\rho
(\omega_{\Gamma_1}(b))= 0$ for all $b \in \B$. Let $\Gamma_0$ denote 
the quotient FODC of $\Gamma_1$ whose associated right ideal 
$\R_{\Gamma_0}$ is generated by $\R_{\Gamma_1}$ and $\B$. Then $\rho$ 
passes to a homomorphism of the quotient FODC $\Gamma_0$ to the 
FODC $\Gamma_2$.}

\mn
{\bf Proof.} If $\Gamma_f$ denotes the universal FODC over $\A$, then any 
other left-covariant FODC $\Gamma$ over $\A$ is isomorphic to the 
quotient FODC $\Gamma_f / \A \omega_{\Gamma_f}(\R_\Gamma)$ (see, for 
instance, [6], Proposition 14.1). This implies that the FODC $\Gamma_0$ is 
isomorphic to the quotient FODC 
$\Gamma_1  / \A \omega_{\Gamma_1}(\R_{\Gamma_0})$. Therefore, 
it is sufficient to prove that $\rho(x \omega_{\Gamma_1}(by))=0$ 
for all $b\in \B$ and $x, y \in \A$. Indeed, using the facts that 
$\rho$ is a bimodule homomorphism and $\omega_{\Gamma_1}(by) = 
S(y_{(1)}) \omega_{\Gamma_1}(b)y_{(2)}$ (see formula (14.3) in [6]) 
we obtain
$$
\qquad\qquad\qquad\qquad\rho(x \omega_{\Gamma_1}(by))= \rho(xS(y_{(1)})) 
\rho(\omega_{\Gamma_1}(b))\rho(y_{(2)}) = 
0.\qquad\qquad\qquad\qquad\qquad\qquad
\hfill\qed$$

\bn
{\bf 3. Preliminaries on the quantum group ${\bf SU_q(2)}$}

\mn
From now let $\A$ be the coordinate $\ast$-algebra $\Oo(SU_q(2))$ of the 
compact quantum group $SU_q(2)$ [14,13,8,6]. 
In what follows we shall assume that $0<q<1$. The generators of 
$\A$ are the four entries $a,b,c,d$ of the fundamental matrix and 
the involution of $\A$ is given by 
\begin{equation}\label{3-31}
a^\ast=d,~ b^\ast=-qc,~ c^\ast=-q^{-1} b,~ d^\ast=a.
\end{equation}

\mn
{\bf 3.1. Star representations of $\Oo(SU_q(2))$}

Suppose that $\pi$ is an arbitrary $\ast$-representation of the $\ast$-algebra 
$\A=\Oo(SU_q(2))$. From Proposition 4.19 in [6] or from the description 
of the irreducible $\ast$-representations in [13] it follows that up the 
unitary equivalence the $\ast$-representation $\pi$ is given by the 
following operator-theoretic model:

Let $v$ and $w$ be unitary operators on Hilbert spaces $\G$ and $\Hh_0$, 
respectively. Put $\Hh={\bigoplus\limits^\infty_{n=0}} \Hh_n$, where 
$\Hh_n=\Hh_0$ for all $n\in\bbbn_0$. For $\eta\in\Hh_0$, let $\eta_n$ denote 
the vector of $\Hh$ which has $\eta$ as its $n$-th component and zero 
otherwise. The $\ast$-representation $\pi$ acts on the direct sum 
Hilbert space $\G\oplus \Hh$ and it is determined by the formulas 
\begin{gather}\label{3-32}
\pi(a)=v,~ \pi(d)=v^\ast, ~ \pi(b)=\pi(c)=0 ~\mbox{ on }~ \G,\\
\label{3-33}
\pi(a)\eta_n=\lambda_n\eta_{n-1},~\pi(d) \eta_n = \lambda_{n+1}\eta_{n+1}, \pi(c) \eta_n=q^nw\eta_n,~ \pi(b)\eta_n=-q\pi(c)^\ast\eta_n=-q^{n+1} 
w^\ast\eta_n
\end{gather}
for $\eta\in\Hh_0$ and $n\in\bbbn_0$, where we have set 
$\eta_{-1}:=0$ and 
\begin{equation*}
\lambda_n:=(1-q^{2n})^{1/2}, n\in\bbbn_0.
\end{equation*}
Note that the $\ast$-representation $\pi$ is parametrized by the 
two unitaries $v$ and $w$.

Let $(\pi,F)$ be a commutator representation of a $\ast$-FODC $\Gamma$ of 
$\A=\Oo(SU_q(2))$. Then $(\pi,F)$ is called 
{\it admissible} if there exist linear subspaces $\E\subseteq\ker \pi(c)$ 
and $\D\subseteq\kere~\pi(a)$ such that 
$\pi(d)\E\subseteq\E, \pi(a)\E\subseteq\E, \pi(b)\D_0\subseteq\D_0, 
\pi(d)\D_0\subseteq\D_0$ and the domain 
$\D_F:=\E+\Lin\{\pi(d^n)\D_0;n\in\bbbn_0\}$ is contained in $\D(F)$ and 
is a core for the (unbounded symmetric) operator $F$. The dense linear 
subspace  
$\D_F$ is then invariant under all operators 
$\pi(x), x\in\A$, and so $\D_F$ may be taken as the domain used in the 
definition of the commutator representation $(\pi,F)$ given in 
the preceding section. 
The admissibility of a representation $(\pi,F)$ is a technical condition 
which is essentially used in the proofs of the main results in Section VI. 
This condition becomes rather natural if it is considered in terms of the 
above model for the $\ast$-representation $\pi$. Clearly, we have 
$\G=\kere~\pi(c), \Hh_0=\kere~\pi(a)$ and $\Hh_n=\pi(d^n)\Hh_0$. 
That $(\pi,F)$ is admissible means that there are linear subspaces 
$\E$ of $\G$ and $\D_0$ of $\Hh_0$ such that $v\E=\E, w \D_0=\D_0$ 
and $\D_F:=\E\oplus \D$ is a core for the operator $F$, where 
$\D :=\Lin\{\D_n;n\in\bbbn_0\}$ and $\D_n:=\{\eta_n;\eta\in\D_0\}$ is 
the $n$-th shift of the domain $\D_0$.

\mn
{\bf 3.2 The $3D$-calculus}
 
The $3D$-calculus was introduced by Woronowicz [14]. A short approach 
was given in [12]. Apart from Section 6 below, we use 
only the following facts concerning the $3D$-calculus.

The right ideal of the $3D$-calculus is generated by the following 
six elements (see [14], (2.27), or [6], (14.23)):
\begin{equation}\label{3-36}
b^2,~ c^2,~  bc,~ (a-1)b,~ (a-1)c,~ q^2a+d - (q^2+1).
\end{equation}
The three 1-forms $\omega_0:=\omega_\Gamma (b), \omega_2:=\omega_\Gamma(c), 
\omega_1:=\omega_\Gamma(a)$ form a basis of the vector space 
${_{inv}\Gamma}$ and the bimodule structure of $\Gamma$ is determined by 
the following commutation relations (see [14], p.135, or [6], p.499):
\begin{align}\label{70}
q \omega_j a = a \omega_j, \omega_j b = q b \omega_j, q \omega_j c = 
c\omega_j, \omega_j d = q d \omega_j\mbox{ for } j = 0,2,\\
\label{71}
q^2\omega_1 a=a\omega_1, \omega_1 b=q^2 b\omega_1, q^2\omega_1c=c\omega_1, 
\omega_1 d=q^2 d\omega_1.
\end{align}

\mn
{\bf 3.3. The $4D_\pm$-calculus}

In order to describe the $4D_\pm$-calculus $\Gamma_\pm$, we first restate 
some facts developed in Subsection 14.2.4 in [6]. The quantum 
tangent space $\T_\pm$ of the $4D_\pm$-calculus is expressed therein in 
terms of the generators $E,F,K,K^{-1}$ of the Hopf algebra 
$\hat{U}_q(sl_2)$ (see [6], p.57). 
Let $\varepsilon_+ :=\varepsilon$ and let $\varepsilon_-$ be the character 
of the algebra $\Oo (SU_q(2))$ such that $\varepsilon_- (a) = 
\varepsilon_- (d)=-1$ and $\varepsilon_- (b)=\varepsilon_-(c)=0$. Then 
$\T_\pm$ is spanned by the four linear functionals
\begin{gather*}
X_1=\varepsilon_\pm K^{-2}-\varepsilon, X_2=q^{1/2}\varepsilon_\pm FK^{-1}, 
X_3=q^{-1/2}\varepsilon_\pm EK^{-1},
X_4=\varepsilon_\pm K^2+ \lambda q^{-1} \varepsilon_\pm FE-\varepsilon
\end{gather*}
and we have
\begin{gather}\label{72}
\Delta (X_1)=\varepsilon\otimes X_1+X_1\otimes\varepsilon_\pm K^{-2},\\
\label{73}
\Delta (X_j)=\varepsilon\otimes X_j+X_j\otimes\varepsilon_\pm +
X_1\otimes X_j,
\ \ j=2,3,\\
\label {74}
\Delta (X_4)=\varepsilon\otimes X_4
+\lambda^2q^{-1}X_1\otimes\varepsilon_\pm
FE+X_4\otimes\varepsilon_\pm K^{2}
 +\lambda^2q^{-1/2}(X_2\otimes \varepsilon_\pm EK
+X_3\otimes\varepsilon_\pm  KF).
\end{gather}
\noindent
There is a dual pairing of Hopf algebras $\breve{U}_q(sl_2)$ and 
$\Oo (SU_q(2))$ determined on the generators by the equations
\begin{equation}\label{75}
\langle K,a\rangle =q^{-1/2}, \langle K,d\rangle = q^{1/2}, 
\langle E,c\rangle =\langle F,b\rangle  = 1\mbox{ and zero otherwise }.
\end{equation}
Let $\{\omega_1,\omega_2,\omega_3,\omega_4\}$ be a basis of the vector 
space ${_{inv}}\Gamma$ such that $(X_i,\omega_j)=\delta_{ij}, 
i,j=1,{\dots},4$. We set $\epsilon=+1$ for the $4D_+$-calculus and 
$\epsilon=-1$ for the $4D_-$-calculus. From (\ref{75}) we derive
\begin{gather}\label{76}
\langle X_1, a\rangle = \epsilon q-1,\langle X_1,d\rangle=\epsilon q^{-1}-1,~\langle X_2,~b\rangle = \langle X_3,c\rangle=\epsilon,\\
\label{77}
\langle X_4,a\rangle =\epsilon q^{-1}-1 +\epsilon \lambda^2 q^{-1}, 
\langle X_4,d\rangle =\epsilon q-1.
\end{gather}
The other pairings of $X_i$ with matrix generators $a,b,c,d$ vanish. Since
$\langle X,x\rangle =(X,\omega (x))$ for $X\in \T_\pm$ and $X\in\A$ (see [6], 
formula (14.9)), we 
therefore obtain \begin{gather}\label{78}
\omega_\Gamma (a)=(\epsilon q-1)\omega_1+(\epsilon q^{-1} - 1+\epsilon\lambda^2q^{-1})\omega_4,~\omega_\Gamma(b)=\epsilon\omega_2,\\
\label{79}
\omega_\Gamma (d)=(\epsilon q^{-1}-1)\omega_1 + (\epsilon q -1)\omega_4, ~ \omega_\Gamma (c)=\epsilon\omega_3.
\end{gather}
From the general theory of left-covariant FODC we recall that
\begin{equation}\label{80}
\omega_i x=\sum\nolimits_j x_{(1)} f^i_j(x_{(2)})\omega_j, x\in\A,
\end{equation}
where the functionals $f^j_i$ are determined by the equation 
$\Delta (X_i)=\varepsilon \otimes X_i+\sum_j X_j\otimes f^j_i$. From the 
formulas (\ref{72})--(\ref{74}) we therefore read off that the non-zero 
functionals $f^i_j$ are
\begin{align*}
&f^1_1=\varepsilon_\pm K^{-2}, f^1_2=X_2, f^1_3=X_3, f^1_4 
=\lambda^2 q^{-1} \varepsilon_\pm FE,f^2_2 = f^3_3=\varepsilon_\pm,\\
&f^2_4=\lambda^2q^{-1/2}\varepsilon_\pm EK, f^3_4=\lambda^2 q^{-1/2} \varepsilon_\pm KF, f^4_4 = \varepsilon_\pm k^2.
\end{align*}
Inserting these functionals into (\ref{78})--(\ref{79}) and using the dual pairing 
(\ref{75}) we obtain the following list of commutation relations which 
described the bimodule structure of the $4D_\pm$-calculus:
\begin{align*}
\omega_1a&=\epsilon qa\omega_1+\epsilon b\omega_3+\epsilon\lambda^2q^{-1} 
a\omega_4,~ \omega_1b=\epsilon q^{-1}b\omega_1+\epsilon a\omega_2,\\
\omega_1c&=\epsilon q c\omega_1+\epsilon d\omega_3+
\epsilon\lambda^2q^{-1} c\omega_4,~ \omega_1 d=\epsilon q^{-1}d\omega_1 
+\epsilon c\omega_2,\\
\omega_2a&=\epsilon a \omega_2+\epsilon \lambda^2q^{-1} b\omega_4, 
~\omega_2 b=\epsilon b\omega_2,~\omega_2c=\epsilon c \omega_2+
\epsilon \lambda^2q^{-1} d\omega_4,~ \omega_2 d=\epsilon d\omega_2,\\
\omega_3a&=\epsilon a \omega_3,~\omega_3 b=\epsilon b\omega_3+
\epsilon\lambda^2 q^{-1} a\omega_4,~\omega_3c=\epsilon c \omega_3,
\omega_3 d=\epsilon d\omega_3+\epsilon\lambda^2 q^{-1} c\omega_4,\\
\omega_4a&=\epsilon q^{-1}a\omega_4,~\omega_4 b=\epsilon qb\omega_4,
\omega_4c=\epsilon q^{-1}c\omega_4,~\omega_4 d=\epsilon qd\omega_4.
\end{align*}
Further, we note that
\begin{equation}\label{89}
\omega^\ast_1=-\omega_1, \omega^\ast_2=-\omega_3, \omega^\ast_4=-\omega_4
\end{equation}
and that the right ideal $\R_{\Gamma_\pm}$ admits the 
following nine generators (see [15], p.132, or [6], p.504):
\begin{equation}\label{3-37}
b^2,~ c^2,~ b(a{-}d),~ c(a{-}d),~ a^2{+}q^2d^2{-}(1{+}q^2)(ad{+}q^{-1}bc),
~ z_\pm b,~ z_\pm c,~ z_\pm(a{-}d),~ z_\pm(q^2a{+}d{-}(q^2{+}1)),
\end{equation}
where 
$$
z_\pm:=q^2a+d-\epsilon(q^3+q^{-1}).
$$

Let $\R_{\pm,3}$ denote the right ideal of $\A$ generated by 
$\R_{\Gamma_\pm}$ and the single element $a+\epsilon qd$. Since 
$S(x)^\ast\in\R_{\pm,3}$ for $x\in\R_{\pm,3}$, there exists a 
left-covariant $\ast$-FODC $\Gamma_{\pm,3}$, of $\A$ which is a 
quotient of the $4D_\pm$-calculus $\Gamma_\pm$ and has the 
associated right ideal 
$\R_{\pm,3}$ (by Proposition 14.6 in [6]). Since 
$\langle X_1, a+\epsilon qd \rangle=0$ by 
(\ref{76})--(\ref{77}), it is not difficult to check that the quantum 
tangent space of $\Gamma_{\pm,3}$ is spanned by the three functionals 
$X_1,X_2,X_3$. Hence 
the FODC $\Gamma_{\pm,3}$ has dimension 3. If we set $\omega_4=0$ in the 
above formulas for the $4D_{\pm}$-calculus $\Gamma_\pm$, then we obtain 
the corresponding formulas for the FODC $\Gamma_{\pm,3}$. In particular, 
we see that the $\ast$-FODC $\Gamma_{+,3}$ gives the classical first 
order differential calculus on $SU(2)$ in the limit $q\rightarrow 1$. 
Thus the $\ast$-FODC $\Gamma_{\pm,3}$ on $SU_q(2)$ seems to be of 
interest in itself.

\mn
{\bf Remark 1.} As noted in [4], the $4D_\pm$-calculus $\Gamma_\pm$ is an 
irreducible bicovariant FODC, because it is derived from the fundamental 
corepresentation of $SL_q(2)$ which is irreducible. By definition this 
means that $\Gamma_\pm$ has no non-trivial {\it bicovariant} quotient FODC. 
However, as we have seen, $\Gamma_{\pm,3}$ is a non-trivial 
{\it left-covariant} quotient FODC of $\Gamma_\pm$.

\mn
{\bf 3.4 The 2-dimensional calculus on the quantum sphere $S^2_q$}

The considerations of this subsection are only needed in Section 6 below. 

Let $\Oo(S^2_q)$ denote the unital $\ast$-subalgebra of $\A=\Oo(SU_q(2))$ 
generated by the elements
$$
x_+ := ba, ~ x_- := cd, ~y_0 := bc.
$$
In order to shorten some formulas it is ocassionally convenient to replace 
$y_0$ by the element
\begin{equation*}
x_0 := \lambda_+ bc + 1 = \lambda_+ y_0 + 1.
\end{equation*}
From the formulas for the comultiplication of the generators $a,b,c,d$ 
in $\Oo(SU_q(2))$ we get
\begin{gather}\label{p1}
\Delta(x_+) = a^2 \otimes x_+ + q^{-1}b^2 \otimes x_- 
+\lambda_+ ba \otimes y_0 +ba \otimes 1,\\
\label{p2}
\Delta(x_-) = qc^2 \otimes x_+ + d^2 \otimes x_- 
+\lambda_+ cd \otimes y_0 + cd \otimes 1,\\
\label{p3}
\Delta(y_0) = ac \otimes x_+ + db \otimes x_- 
+x_0 \otimes y_0 +bc \otimes 1.
\end{gather}
That is, we have $\Delta (\Oo(S^2_q)) \subseteq \A \otimes \Oo(S^2_q)$ and 
hence $\Oo(S^2_q)$ is a left quantum space of $\A$ (that is, a 
left $\A$-comodule algebra) with coaction given by the restriction of 
the comultiplication. 
It is well-known that that $\Oo(S^2_q)$ is the coordinate algebra of the 
Podles' quantum 2-sphere $S^2_{qc}$ in the case $c=0$ [9]. 
(Note that the quantum 2-spheres in [9], [10] and [1] are right quantum spaces, 
while we consider the corresponding left quantum spaces here.) The 
generators $x_+,x_-,y_0$ satisfy the relations 
\begin{gather}\label{p3a}
x_+x_- -q^2 x_-x_+ =(q^2-1)y_0^2 ,~ x_+x_- -q^4 x_-x_+ =(1-q^2)q y_0, \\
\label{p3b}
x_+y_0 = q^2 x_+y_0, q^2 x_-y_0 = y_0x_-.
\end{gather}
In fact, the algebra $\Oo(S^2_q)$ can be 
also characterized as the abstract unital algebra with 
generators $x_+, x_-, y_0$ and definining relations 
(\ref{p3a}) and (\ref{p3b}). Since $q$ is real, 
the algebra $\Oo(S^2_q)$  is a $\ast$-algebra with 
involution determined by $(x_+)^\ast = x_-$ and $(y_0)^\ast = y_0$.

As shown by P. Podles [10], the quantum space $S^2_q$ carries a unique 
2-dimensional $\ast$-FODC. For the application given in Section 6 
it is crucial that this calculus is induced from the $3D$-calculus of 
the quantum group $SU_q(2)$. This fact has been known to the author 
since several years (in fact, since the writing of [1])  
and also to others (S. Majid, P. Podles). Since 
I could not find this result in the literature, 
we shall derive it in this subsection. In order to do so we first repeat 
some more facts on the $3D$-calculus from Subsection 14.1.3 in [6].  
Let $(\Gamma,\dd)$ be the $3D$-calculus on $\Oo(SU_q(2))$ and let 
$\Gamma_2$ denote the induced $\ast$-FODC of $\Gamma$ on the 
$\ast$-subalgebra $\Oo(S^2_q)$.

The quantum tangent space $\T$ of the $3D$-calculus has the three basis 
elements
\begin{equation*}
X_0:=q^{-1/2} FK,\ \ \ X_2:=q^{1/2} EK,\ \ \ X_1:=(1- q^{-2})^
{-1}(\varepsilon-K^{4})
\end{equation*}
satisfying
\begin{equation}\label{p4}
\Delta X_j=\varepsilon\otimes X_j+X_j\otimes K^{2}\ {\rm for}\ j=0,2\ \
{\rm and}\ \ \Delta X_1=\varepsilon\otimes X_1+X_1\otimes K^{4}. 
\end{equation}
The basis $\{\omega_0,\omega_1,\omega_2\}$ of the vector space 
${_{\rm inv}\Gamma}$ is dual to the basis $\{X_0,X_1,X_2\}$ of $\T$. 
Therefore, by the general theory of left-covariant 
differential calculi we have
\begin{equation}\label{p5}
\dd x ={\sum^2_{j=0}} x_{(1)} \langle X_j, x_{(2)} \rangle \omega_j,
\ \ \ x\in \A.
\end{equation}
From (\ref{p4}) and the relations 
$\langle X_0,b \rangle =\langle X_2,c \rangle = \langle X_1,a \rangle =1$ and 
$\langle X_0,a \rangle = \langle X_0,c \rangle= 
\langle X_2,a \rangle= \langle X_2,b \rangle= \langle X_1,b \rangle= 
\langle X_1,c \rangle = 0$ by (\ref{75}), we obtain
\begin{gather*}
\langle X_0,x_+ \rangle =\langle X_2,x_- \rangle = q^{-1},\\
\langle X_0,x_- \rangle = \langle X_0,y_0 \rangle = 
\langle X_2,x_+ \rangle = \langle X_2,y_0 \rangle =
\langle X_1,x_+ \rangle= \langle X_1,x_- \rangle = 
\langle X_1,y_0 \rangle = 0.
\end{gather*} 
Inserting these facts and equations (\ref{p1})--(\ref{p3}) into (\ref{p5}) 
we get
\begin{gather}\label{p6}
\dd x_+ = q^{-1}a^2 \omega_0 + b^2 \omega_2,~
\dd x_- = c^2 \omega_0 + q d^2 \omega_2,~
\dd y_0 = ca \omega_0 + bd \omega_2.
\end{gather}
Some lengthy but straightforward 
computations using the formulas (\ref{p6}),(\ref{70}) and (\ref{71}) yield the following commutation relations for the 
FODC $\Gamma_2$ on $\Oo(S^2_q)$:
\begin{gather*}
\dd x_+ x_+ = x_+ \dd x_+ -q^{-1}\lambda x_+^2 \dd x_0 + 
q\lambda x_+x_0 \dd x_+,\\
\dd x_+ x_- = q^2 x_- \dd x_+ + q\lambda x_+x_- \dd x_0 -
q^{-1}\lambda x_+(x_0{-}1) \dd x_-,\\
\dd x_+ x_0 = x_0 \dd x_+ + q\lambda x_+(x_0{+}q^{-2}) \dd x_0 
-q^{-1} \lambda \lambda_+^2 x_+^2 \dd x_-,\\
\dd x_- x_+ = q^{-2} x_+ \dd x_- -q^{-1} \lambda x_-x_+ \dd x_0
-q \lambda x_-(x_0{-}1) \dd x_+,\\
\dd x_- x_- = x_- \dd x_- + q \lambda x_-^2 \dd x_0 
- q^{-1} \lambda x_- x_0 \dd x_- ,\\
\dd x_- x_0 =  x_0 \dd x_- -q^{-1} \lambda x_-(x_0{+}q^2) \dd x_0
+q\lambda \lambda_+^2 x_-^2 \dd x_+,\\
\dd x_0 x_+ =q^{-2} x_+ \dd x_0 +q^{-1}\lambda x_+(x_0{+}q^{-2}) \dd x_0 
- q^{-1}\lambda (x_0{-}1) \dd x_+ -q^3 \lambda \lambda_+^2 x_+^2 \dd x_-,\\
\dd x_0 x_- = q^2 x_- \dd x_0  +q\lambda(x_0{-}1) \dd x_- 
-q\lambda x_-(x_0{+}q^2) \dd x_0 +q^3 \lambda \lambda_+^2 x_-^2 \dd x_+,\\
\dd x_0 x_0 = x_0 \dd x_0 - q^{-1} \lambda \lambda_+^2 (x_0{-}1)x_+ \dd x_- +
q\lambda x_0(x_0{-}1) \dd x_0.
\end{gather*}
Moreover, from (\ref{p6}) it follows also that
\begin{equation}\label{p7}
x_+ \dd x_- + q^2 x_- \dd x_+ - q x_0 \dd y_0 = 0.
\end{equation}
\mn
{\bf Lemma 4.} {\it The FODC $\Gamma_2$ is the left counter-part of the 
unique 2-dimensional left-covariant 
$\ast$-FODC on $\Oo(S^2_q)$ characterized in [10].} 

\mn
{\bf Proof.} The assertion will follow from the first statement 
of the main theorem in [10]. Hence it suffices to check that 
$\Gamma_2$ fulfills the assumptions made there. Being induced from 
the left-covariant $\ast$-FODC $\Gamma$ on $\Oo(SU_q(2))$, $\Gamma_2$
is obviously a left-covariant $\ast$-FODC on $\Oo(S^2_q)$. The 
above commutation rules show that the differentials 
$\dd x_+, \dd x_-, \dd x_0$ generate $\Gamma_2$ as a left 
$\Oo(S^2_q)$--module. Thus it remains to verify assumption 7) in [10], Section 1. In 
the present context this condition means that for arbitrary elements 
$z_+, z_-, z_0 \in \Oo(S^2_q)$ an equation 
\begin{equation}\label{p15}
z_+ \dd x_+ + z_- \dd x_- + z_0 \dd y_0 = 0
\end{equation}
in $\Gamma_2$ is valid if and only if there is an element $z \in \Oo(S^2_q)$  
such that 
\begin{equation}\label{p14}
 z_+ = q^2 zx_- ,~ z_- = zx_+ ,~ z_0 = -q zx_0.
\end{equation}
Clearly, (\ref{p14}) implies (\ref{p15}) because of the relation 
(\ref{p7}). (By modifying the uniqueness proofs given in 
[10] or [1] the proof of the converse direction can be avoided in 
the present case, but we prefer to carry out it here.) 
Conversely, suppose now that (\ref{p15}) holds. 
Inserting this into (\ref{p6}) and comparing the coeffients of 
$\omega_0$ and $\omega_2$ we obtain 
\begin{gather}\label{p8}
q^{-1} z_+ a^2 + z_- c^2 + z_0 ca = 0,\\
\label{p9}
z_+ b^2 + q z_- d^2 + z_0 bd = 0.
\end{gather}
The equations (\ref{p9})$q^2ac$--(\ref{p8})$db$, 
(\ref{p9})$a^2$--(\ref{p8})$q^{-3}b^2$ and 
(\ref{p8})$d^2$--(\ref{98})$q^3c^2$ can be written as 
\begin{gather}\label{p10}
q^2z_-x_- - z_+x_+ + q\lambda z_0y_0 = 0,\\
\label{p11}
z_-(q^{-1}\lambda_+ y_0 + q) + z_0x_+ = 0,\\
\label{p12}
z_+(q\lambda_+ y_0 + q^{-1}) + z_0x_- = 0,
\end{gather}
respectively. Define now an element $z \in \Oo(S^2_q))$ by
\begin{equation}\label{p13}
z:= - \lambda_+^2 z_-x_- - z_0(q\lambda y_0 + q^{-1}) =
 - q^{-2}\lambda_+^2 z_+x_+ - z_0(q^{-3}\lambda_+y_0 + q^{-1}),
\end{equation}
where the second equality follows from (\ref{p10}). From the 
algebra relations (\ref{p3a})--(\ref{p3b}) and the 
formulas (\ref{p11}) and 
(\ref{p12}) it then follows that the relations (\ref{p14}) are fulfilled. 
For instance, let us explain how to get the first equality 
of (\ref{p14}). First we multiply the second 
representation of $z$ in (\ref{p13}) by $x_-$ from the right, then 
we symplify the terms by means of the algebra 
relations $x_+x_- = q^2y_0^2+qy_0$ 
and $y_0x_-=q^2x_-y_0$ and finally we insert the expression of 
$z_0x_-$ from (\ref{p12}). This in turn yields 
the desired relation $zx_- = q ^{-2}z_+$. The second and third equalities 
in (\ref{p14}) can be derived in similar manner from 
the first expression of $z$ in (\ref{p13}) and formula (\ref{p11}).
\hfill\qed

\bn
{\bf 4. Main results}

\mn
The first main theorem describes all possible admissible 
commutator representations of the $3D$-calculus on $SU_q(2)$. In order to 
formulate this result some further preliminaries are needed.

Let $\pi$ be a $\ast$-representation of $\A=\Oo(SU_q(2))$ as described 
by the model in the preceding section. Suppose that are a linear 
operator $T$ and a symmetric linear operator $R$ on $\Hh_0$ and a dense 
linear subspace $\D_0\subseteq\D(T) \cap\D(T^\ast)\cap\D(R)$ of $\Hh_0$ 
such that $w\D_0=\D_0$,
\begin{gather}\label{4-41}
w Tw^\ast\eta=q T\eta ~\mbox { for }~ \eta\in \D_0,\\
\label{4-42}
w^2Rw^{\ast 2} \eta+q^2 R\eta=(1+q^2)wRw^\ast \eta~ \mbox { for } ~\eta \in\D_0. 
\end{gather}
From the assumptions $w\D_0=\D_0$ and (\ref{4-41}) one easily derives 
that \begin{equation}
wT\eta = q Tw\eta,~ wT^\ast \eta = q T^\ast w \eta, 
~ Tw^\ast \eta = q w^\ast T\eta,~ T^\ast w^\ast\eta  = q w^\ast T^\ast \eta ~ {\rm for}~
\eta \in \D_0.
\end{equation}
Further, assume that there are a symmetric linear operator $Q$ on $\G$ 
and a dense linear subspace $\E$ of $\G$ such that $v\D_0=\D_0$ and 
\begin{equation}\label{4-43}
v^2 Qv^{\ast 2}\eta+q^2 Q\eta=(1+q^2)vQv^\ast\eta ~\mbox { for }~ 
\eta\in \D_0.
\end{equation}
Let $F$ be a linear operator on the Hilbert space $\G\oplus \Hh$ which 
has the dense linear subspace 
$\D_F:=\E\oplus \Lin\{\eta_n:\eta\in\D_0, n\in\bbbn_0\}$
as a core and is defined by
\begin{gather}\label{4-44}
F \eta_n=\lambda_n T\eta_{n-1}+w^nRw^{\ast n} \eta_n+\lambda_{n+1} T^\ast 
\eta_{n+1},~ \eta\in\D_0,\\
\label{4-45}
F\eta = Q\eta,~ \eta\in\E.
\end{gather}
Clearly, $F$ is a symmetric operator.

\mn
{\bf Theorem 1.} {\it Under the above assumptions, the pair $(\pi,F)$ is an 
admissible commutator representation of the $3D$-calculus on $SU_q(2)$. 
Up to unitary equivalence any admissible commutator representation of 
the $3D$-calculus is of this form.}
\mn

We shall see in the next section that by appropriate choice of the above 
operators $T$ and $R$ one obtains a faithful admissible commutator 
representation of the $3D$-calculus. In contrast to this the 
$4D_\pm$-calculus has no faithful admissible commutator representation.

\bn
{\bf Theorem 2.} {\it Let $(\pi,F)$ be an admissible commutator 
representation of the $4D_\pm$-calculus $\Gamma_\pm$. 
Then the corresponding $\ast$-FODC homomorphism 
$\rho:\Gamma_\pm \rightarrow\Gamma_{\pi,F}$ passes to a homomorphism of 
the quotient $\ast$-FODC $\Gamma_{\pm, 3}$ to 
$\Gamma_{\pi,F}$, so $(\pi,F)$ becomes a commutator representation 
of $\Gamma_{\pm, 3}$.}

\bn
The next theorem shows in particular that none of the three calculi can 
be given by a spectral triple in the sense of A. Connes.

\bn
{\bf Theorem 3.} {\it If $(\pi, F)$ is a commutator representation of the 
$3D$-calculus or the $4D_\pm$-calculus such that 
all operators $\dd_{\pi,F}(x), x\in\A$, are bounded, then we have 
$\dd_{\pi, F}(x)=0$ for all $x\in\A$.}

\bn
{\bf 5. Commutator representations of the $3D$-calculus on $SU_q(2)$}

\mn
In this section we investigate admissible commutator representations 
$(\pi, F)$ of the $3D$-calculus $\Gamma$ more in detail. 
Throughout this section we retain the notation of 
Sections 3 and 4 and suppose that $\pi$ is a 
$\ast$-representation of $\A=\Oo(SU_q(2))$ such that $\G=\{0\}$. If 
not specified otherwise all operator equations containing the 
operators $R,T$ and $F$ are meant on the domains $\D_0$ and 
$\D_F=\Lin\{\eta_n;\eta\in\D_0\}$, respectively. 
Further, we will denote an operator on $\Hh_0$ and the corresponding 
diagonal operator on $\Hh={\bigoplus\nolimits_n} \Hh_n$ by the same symbol.

Let us first look at the operator relation (\ref{4-42}). It can be rewritten in the form
\begin{equation}\label{6.1}
w(wRw^\ast-R)w^\ast=q^2(wRw^\ast-R).
\end{equation}
Therefore, if $R^\prime$ and $R^{\prime\prime}$ satisfies the operator 
equations 
\begin{equation}\label{6.2}
wR^\prime w^\ast=q^2R^\prime~{\rm and}~ wR^{\prime\prime}w^\ast=
R^{\prime\prime},
\end{equation}
respectively, then $R:=R^\prime+R^{\prime\prime}$ 
is a solution of equation (\ref{4-42}). Conversely, suppose that $R$ is a 
solution of (\ref{4-42}) and put
$$
R^\prime:=(1+q^2)^{-1}(R-wRw^\ast)~ {\rm  and}~R^{\prime\prime}=(1+q^2)^{-1}(q^2R+wRw^\ast).
$$
Then $R^\prime$ and $R^{\prime\prime}$ satisfy the 
equations (\ref{6.2}) and we have $R:=R^\prime+R^{\prime\prime}$. 
From this decomposition and Lemma 6
 below it follows in particular that the only {\it bounded} solutions of (\ref{4-42}) are the bounded operators commuting with $w$. Further, we obtain that 
\begin{equation}\label{6.4}
R_n\eta_n\equiv w^n R w^{\ast n} \eta_n = q^{2n} R^\prime \eta_n + 
R^{\prime\prime} \eta_n, ~\eta\in \D_0.
\end{equation}

If $(\pi,F)$ is a commutator representation of $\Gamma$ with $\ast$-FODC 
homomorphism $\rho:\Gamma\rightarrow \Gamma_{\pi,F}$, we let 
$\Omega_j=\rho(-{\rm i}\omega_j)$ denote the image of 
the left-invariant 1-form $-{\rm i}\omega_j, j=0,1,2$. Then we have
$$
\Omega_0=\Omega_{\pi,F}  (b),~ \Omega_2=\Omega_{\pi,F}(c),~
\Omega_1=\Omega_{\pi,F}(a)=-q^{-2} \Omega_{\pi,}(d).
$$

The next theorem gives a reformulation of admissible commutator 
representations in terms of the representation $\pi$. It shows that 
the operators $\Omega_0,\Omega_2,\Omega_1$ can be nicely expressed in 
terms of the operators $T$ and $R^\prime$. Note that $\pi (c)^{-1}$ is 
a well-defined bounded operator mapping $\D_F$ into itself, because we 
assumed that $\G=\{0\}$.

\bn
{\bf Theorem 4.} {\it Suppose that $(\pi,F)$ is an admissible commutator 
representation of the $3D$-calculus and let $F$ be of the form 
(\ref{4-44}) with 
operators $T$ and $R=R^\prime + R^{\prime\prime}$ satisfying (\ref{4-41}) 
and (\ref{4-42}), respectively. Then we have
\begin{gather}\label{6.4b}
\pi (c)T=qT\pi(c)~ and~\pi(c)R^\prime=q^2R^\prime\pi(c),\\
\label{6.5}
R^{\prime\prime} \pi(x)=\pi(x)R^{\prime\prime}~for~all~x\in\A,\\
\label{6.6}
F\eta_n= T\pi(a)\eta_n+\pi(c)^nR^\prime\pi(c)^{-n}\eta_n+R^{\prime\prime} 
\eta_n+T^\ast\pi(d)\eta_n,~\eta\in\D_0,\\
\label{6.7}
\Omega_0=\lambda \pi(b)T,~~\Omega_2=- \lambda\pi(c)T^\ast, 
~~\Omega_1= q^{-2} \lambda \pi(bc)R^\prime.
\end{gather}
Conversely, if $R^\prime$ and $R^{\prime\prime}$ are symmetric linear 
operators and $T$ is a linear operator defined on common dense linear 
subspace $\D_0\subseteq \D(T) \cap\D(T^\ast)\cap\D(R)$ 
of the Hilbert space $\Hh_0$ such that (\ref{6.4b}) and (\ref{6.5}) are 
valid, then the pair $(\pi,F)$ with $F$ defined by (\ref{6.6}) is an 
admissible commutator representation of the $3D$-calculus.}

\bn
{\bf Proof.} Most of the assertions are only reformulations of the 
conditions occuring in Section 4.  Therefore we do not carry out all 
details of proof. For instance, (\ref{4-41}) and (\ref{4-42}) are 
equivalent to the equations (\ref{6.4b}) and (\ref{6.5}). Since 
$R^{\prime\prime}=wR^{\prime\prime} w^\ast$ as noted above, 
$R^{\prime\prime}$ commutes with $\pi(b)$ and $\pi(c)$ and hence with 
all representation operators $\pi(x),x\in\A$. Formula (\ref{6.6}) follows 
from (\ref{4-44}) and (\ref{6.4}).

As a sample, we prove the formula for $\Omega_0$ and compute 
\begin{align*}
\Omega_0\eta_n=~&(\pi(d)F\pi(b)-q^{-1}\pi(b) F\pi(d))\eta_n\\
=~&-q^{n+1}\pi(d)(\lambda_nTw^\ast\eta_{n-1}+w^nRw^{\ast n+1}\eta_n+
\lambda_{n+1}T^\ast w^\ast\eta_{n+1})\\
&~-q^{-1}\pi(b)\lambda_{n+1}(\lambda_{n+1}T\eta_n+w^{n+1}Rw^{\ast n+1} 
\eta_{n+1}+\lambda_{n+2} T^\ast \eta_{n+2})\\
=~&-q^{n+1}(\lambda^2_n Tw^\ast \eta_n+\lambda_{n+1} w^n R w^{\ast n+1} 
\eta_{n+1} + \lambda_{n+1}\lambda_{n+2} T^\ast w^\ast\eta_{n+2})\\
&+q^n \lambda^2_{n+1} w^\ast T\eta_n+q^{n+1}
\lambda_{n+1} w^n R w^{\ast n+1}\eta_{n+1}
+q^{n+2} \lambda_{n+1} \lambda_{n+2} w^\ast T^\ast \eta_{n+2}\\
=~&(-q^{n+1} \lambda^2_n+q^{n-1} \lambda^2_{n+1}) Tw^\ast \eta_n\\
=~&q^n(q^{-1}-q) T w^\ast\eta_n= \lambda \pi (b)T \eta_n.
\end{align*}
for $\eta\in\D_0$. In similar manner, one shows that
\begin{align*}
\Omega_2\eta_n=(-q\pi(c)F\pi(a)+\pi(a)F\pi(c))\eta_n=
q^{n-1}(q^{-1}-q)wT^\ast\eta_n,
\end{align*}
\begin{align*}
\Omega_1 \eta_n &=
( \pi (d)F \pi(a) - q^{-1} \pi (b) F \pi (c) - F) \eta_n 
= (R_{n -1} - R_n) \eta_n\\
&= w^{n-1} (R-w R w^\ast) w^{\ast n-1} \eta_n 
= (1 - q^2) w^{n-1} R^\prime w^{\ast n-1} \eta_n = 
(1 - q^2) q^{2n-2} R^\prime \eta_n.
\end{align*}
These relations imply the two other formulas of (\ref{6.7}).
\hfill\qed

\sn
By the preceding, for a given $\ast$-representation $\pi$ of $\A$ 
such that $\G=\{0\}$ the operators $F$ of admissible pairs $(\pi,F)$ are 
parametrized by the three operators $T,R^\prime$ and $R^{\prime\prime}$ on 
the Hilbert space $\Hh_0$ satisfying the relations
\begin{equation}\label{6.8}
wTw^\ast=qT, ~wR^\prime w^\ast=q^2R^\prime~
{\rm and}~ wR^{\prime\prime}w^\ast=R^{\prime\prime}.
\end{equation}

It is now easy to construct admissible pairs $(\pi,F)$. We shall do this 
for the faithful $\ast$-representation $\pi$ of the $\ast$-algebra $\A$ 
given in [14]. In this case $w$ is the 
backward shift on the Hilbert space $\Hh_0=l_2(\bbbz)$. That is, if we 
identify $\Hh$ with $l_2(\bbbn_0 \times\bbbz)$ and denote by 
$\{e_{nk}; n\in\bbbn_0, k\in\bbbz\}$ the standard orthonormal basis 
of $l_2 (\bbbn_0\times\bbbz)$, then the operators $w,\pi(a)$ and 
$\pi(c)$ act as 
\begin{equation}\label{6.9}
w e_{nk}=e_{n,k-1},~ \pi(a)e_{nk}=\lambda_n e_{n-1,k},~
 \pi(c)e_{nk}=q^n e_{n,k-1}.
\end{equation}
Define linear operators $T$ and $R^\prime$ on the domain 
$\D_F:={\rm Lin} \{e_{nk}; n\in\bbbn_0, k\in\bbbz\}$ by
\begin{equation}\label{6.10}
T e_{nk} = q^k e_{n,  k-1},~R^\prime e_{nk} = q^{2k} e_{kn}.
\end{equation}
Let $R^{\prime\prime}$ by a symmetric linear operator on $\D_F$ such 
that $wR^{\prime\prime}w^\ast=R^{\prime\prime}$. The conditions 
(\ref{6.8}) are obviously fulfilled. By (\ref{6.4}), 
(\ref{6.10}) and (\ref{4-44}), 
the action of the operator $F$ on the basis vectors $e_{nk}$ is given by
$$
F e_{nk}=\lambda_n q^k e_{n-1,k-1}+ q^{2n+2k}e_{nk}+\lambda_{n+1} 
q^{k+1} e_{n+1,k+1} +R^{\prime\prime} e_{nk}.
$$
{\it Then the pair $(\pi,F)$ is an admissible commutator representation of 
the $3D$-calculus $\Gamma$ on $SU_q(2)$.} For instance, one may take 
$R^{\prime\prime}$ of the form $R^{\prime\prime}e_{nk}=
\sum_r\alpha_r e_{n,k-r}$, where $(\alpha_r; r\in\bbbz)$ is a real sequence 
such that $\alpha_r=0$ for $|r|\ge r_0$. In this case it is 
straightforward to prove that then the corresponding FODC homomorphism 
$\rho:\Gamma\rightarrow\Gamma_{\pi,F}$ is {\it faithful}. 
(Indeed, using the vector space basis $\{a^nb^mc^r, b^mc^rd^s;m,n,s\in 
\bbbn_0,n\in\bbbn\}$ of $\A$ and the formulas (\ref{6.7}) one verifies that 
any relation $\pi(x_0)\Omega_0+\pi(x_1)\Omega_1+\pi(x_2)\Omega_2=0$ with 
$x_0,x_1,x_2\in\A$ implies that $x_0=x_1=x_2=0$.) 

Note that the operators ${\rm d}_{\pi,F} (x)=[{\rm i}F,\pi (x)],x\in\A$, of 
the FODC $\Gamma_{\pi,F}$ are unbounded. This stems from the fact that 
$T$ and $R^{\prime\prime}$ and hence the basis elements $\Omega_j, j=0,1,2$, 
of the vector space of left-invariant 1-forms of $\Gamma_{\pi,F}$ are 
unbounded operators. The reason are the sequences $(q^k)$ resp. $(q^{2k})$ 
in (\ref{6.10}) as $k\rightarrow -\infty$, so this unboundedness 
is rather well controlled. 

Commutatator representations can be used to construct extensions 
of $\ast$-FODC to larger algebras. We explain this for the $3D$-calculus. It is clear that the set 
$\Ss:= \{b^n c^m; n,m\in \bbbn_0 \}$ is a left Ore subset of the 
algebra $\A$ (that is, for any $(s,x) \in \Ss \times \A$ there exists 
$(t,y) \in \Ss \times \A$ such that $ys=tx$.) Moreover, the algebra $\A$ has no 
zero divisors. Therefore, as it is well-known in ring theory, there 
exists a $\ast$-algebra ${\tilde \A}$ which contains $\A$ as a 
$\ast$-subalgebra such that the elements of $\Ss$ are invertible 
and ${\tilde \A}$ is generated by $\A$ and the inverses of $\Ss$. 
Since the Hilbert space $\G$ is zero, the $\ast$-representation $\pi$ 
of $\A$ extends uniquely to a $\ast$-representation ${\tilde \pi}$ of the 
$\ast$-algebra ${\tilde \A}$. Hence 
$(\Gamma_{{\tilde \pi},F},\dd_{{\tilde \pi},F})$ is a $\ast$-FODC of the 
$\ast$-algebra ${\tilde \A}$. If we take a faithful 
commutator representation $(\pi,F)$ of the $3D$-calculus, the we obtain 
an extension of the $3D$-calculus to the larger $\ast$-algebra ${\tilde \A}$ 
in this manner.  

For the study of harmonic analysis and metric noncommutative geometry 
on $SU_q(2)$ it is more important to work with the direct sum 
$\pi_{\rm reg}$ of $\bbbn_0$ copies of the 
$\ast$-representation $\pi$. That is, we take the Hilbert space 
$\Hh_{\rm reg}=l_2(\bbbn \times\bbbz 
\times\bbbn_0)$ with standard orthonormal basis 
$\{ e_{nkl} ; n,l,\in\bbbn_0, k\in\bbbz\}$ and let the operators 
$\pi_{\rm reg} (x), x\in\A$, and $w$ act on the first two indices as 
stated above. Then $\pi_{\rm reg}$ is just the GNS representation of $\A$ 
associated with the Haar state $h$ of the compact quantum group algebra 
$\A=\Oo (SU_q(2))$. Indeed, if $\varphi_h$ denotes the vector
\begin{equation}\label{6.11a}
\varphi_h :=(1-q^2)^{-1/2} \sum^\infty_{n=0} q^n e_{n0n},
\end{equation}
then it follows at once from the explicit formulas for the 
Haar state [14,13,8,6] 
that 
\begin{gather*}
h(x)=\langle \pi_{\rm reg} (x) \varphi_h, \varphi_h\rangle, ~x\in\A.
\end{gather*}
Let $\alpha$ and $\beta$ be positive reals. Define 
operators $T, R^\prime$ and $R^{\prime\prime}$ 
on the span of basis vectors 
$e_{nkl}$ by
\begin{equation}\label{6.11b}
T e_{nkl} = \alpha (1+q^2)^{1/2} q^k e_{n,k-1,l-1},
~R^\prime e_{nk} = \beta q^2(1+q^2+q^4)^{1/2} q^{2k} e_{knl}, 
~R^{\prime\prime} = 0.
\end{equation}
Let $F_{\rm reg}$ denote the 
corresponding operator given by (\ref{6.6}). Then the pair 
$(\pi_{\rm reg},F_{\rm reg})$ is another admissible 
commutator representation  of the $3D$-calculus. It is  
natural to use the state vector $\varphi_h$ to define a scalar 
product on the 1-forms of the $3D$-calculus $\Gamma$ by 
\begin{equation}\label{6.11c}
\langle \omega, \omega^\prime \rangle :=
\langle\rho (\omega)\varphi_h,\rho(\omega^\prime)\phi_h \rangle, ~~\omega, 
\omega^\prime \in \Gamma,
\end{equation}
where $\rho :\Gamma\rightarrow \Gamma_{\pi_{\rm reg}, F_{\rm reg}}$ is 
the corresponding $\ast$-FODC homomorphism. Using the 
formulas (\ref{6.11a}), (\ref{6.11b}) and (\ref{6.11c}) we compute 
$$
\langle \omega_0, \omega_0 \rangle = 
\langle \omega_2, \omega_2 \rangle = \alpha^2, ~
\langle \omega_1, \omega_1 \rangle = \beta^2
~{\rm and }~\langle \omega_k, \omega_l \rangle = 0 
~{\rm if}~ k \neq l.
$$ 

\bn
{\bf 6. Commutator representations of the 2-dimensional calculus on $S^2_q$}

By Lemma 4 we have shown that the $3D$-calculus $\Gamma$ on the 
quantum group $SU_q(2)$ induces the 2-dimensional calculus $\Gamma_2$ on the 
quantum 2-sphere $S^2_q$. Thus any commutator representation of the $3D$-calculus gives obviously a commutator representation of the 
$\ast$-FODC $\Gamma_2$ on the $\ast$-subalgebra $\Oo(S^2_q)$ 
of $\Oo(SU_q(2)$. In this brief section we shall make this more explicit.

Let $(\pi,F)$ be an admissible commutator representation of the $3D$-calculus 
as described in Section 4, where the $\ast$-representation $\pi$ is 
as in Subsection 3.1. Using the condition $wTw^\ast = qT$ by (\ref{4-41}) 
and the formulas (\ref{3-32}) and (\ref{3-33}) 
we compute the differentials of 
the generators $x_+, x_-, y_0$ and obtain
\begin{gather*}
\dd_{\pi,F}(x_+) = {\rm i}q^{-1}\lambda \pi(a)T \pi(x_+)
- {\rm i} \lambda \pi(b) T^\ast \pi (y_0),~~
\dd_{\pi,F}(x_-) = -{\rm i}q\lambda \pi(d)T \pi(x_-)
+ {\rm i} \lambda \pi(c) T \pi (y_0),\\
\dd_{\pi,F}(y_0) = {\rm i}q^{-1}\lambda \pi(a)T \pi(y_0)
- {\rm i} q \lambda \pi(d) T \pi (y_0) 
= {\rm i}\lambda \pi(c)T \pi(x_+)
- {\rm i} \lambda \pi(b) T^\ast \pi (x_-).
\end{gather*}
These formulas describe the corresponding commutator representation of the 
$\ast$-FODC $\Gamma_2$ of $\Oo(S^2_q)$. In particular we see that the 
operator $R$ does not occur in these formulas and that 
$\dd_{\pi,F} (x) = 0$ on the subspace $\G$ for all 
$x \in \Oo(S^2_q)$. 

Conversely, let $\pi$ be a $\ast$-representation of $\Oo(SU_q(2))$ 
on the Hilbert space $\Hh = \oplus_n \Hh_n$ as 
in Subsection 3.1 and let $T$ be a linear operator on the 
Hilbert space $\Hh_0$. If there exists a 
dense linear subspace $\D_0\subseteq\D(T) \cap\D(T^\ast)$ of $\Hh_0$ 
such that $w\D_0=\D_0$ and $wTw^\ast \eta  = qT \eta$ for $\eta \in \D_0$, 
then the above formulas describe a commutator representation of the $\ast$-FODC 
$\Gamma_2$ of $\Oo(S^2_q)$. Examples can be constructed similarly as in the 
case of the $3D$-calculus. 
 
Using the $\ast$-FODC $(\Gamma_{{\tilde \pi},F},\dd_{{\tilde \pi},F})$ 
of the Ore extension ${\tilde \A}$ of the $\ast$-algebra $\A$ defined in the 
preceding section, the operators $T$ and $T^\ast$ 
can be expressed by the formulas
\begin{equation}\label{p30}
T = {\rm i}\lambda^{-1} b \dd_{{\tilde \pi},F}(db^{-1}) ~{\rm and}~ 
T^\ast = -{\rm i}\lambda^{-1} c \dd_{{\tilde \pi},F}(ac^{-1}). 
\end{equation}

\mn
{\bf Remark 2.} Let $z:=ac^{-1}$. In the $\ast$-algebra ${\tilde \A}$ 
we then have $z^\ast = -db^{-1}$ and 
\begin{equation}\label{p31}
z^\ast z - q^2 z z^\ast = q^2-1.
\end{equation}
The $\ast$-subalgebra $\Z$ of ${\tilde \A}$ generated 
by the element $z:=ac^{-1}$ is just the abstract 
$\ast$-algebra with a single generator $z$ 
and defining relation (\ref{p31}). It is well-known that this 
$\ast$-algebra $\Z$ has a $\ast$-FODC with 
commutation relations 
\begin{equation*}
\dd z {\cdot} z = q^2 z \dd z,~ \dd z {\cdot} z^\ast = q^{-2} z^\ast \dd z,~
\dd z^\ast {\cdot} z = q^2 z \dd z^\ast,~ 
\dd z^\ast {\cdot} z^\ast = q^{-2} z^\ast \dd z^\ast.
\end{equation*}
These relations can be found (for instance) in [11]. 
As a byproduct of the preceding consideration we obtain a commutator 
representation $(\pi,F)$ of this $\ast$-FODC, where $\pi$ denotes the 
restriction to $\Z$ of the above $\ast$-representation ${\tilde \pi}$ 
of ${\tilde \A}$ 
and $F$ is the operator given by (\ref{4-44}) with $R=0$ and $wTw^\ast = qT$. 
That is, we have
\begin{gather*}
\pi(z) \eta_n = q^{-n} w^\ast \lambda_n \eta_{n-1}, ~
\pi(z^\ast) \eta_n = q^{-n-1} w \lambda_{n+1} \eta_{n+1},\\
F \eta_n = \lambda_n T \eta_{n-1} + \lambda_{n+1} T^\ast \eta_{n+1}.
\end{gather*}
From formula (\ref{p30}) we see that then the operators 
of the differentials 
$\dd_{\pi,F}(z)$ and $\dd_{\pi,F}(z^\ast)$ act as
\begin{equation*}
\dd_{\pi,F}(z) \eta_n = {\rm i} \lambda q^{-n} w^\ast T^\ast \eta_n,~~
\dd_{\pi,F}(z^\ast) \eta_n = -{\rm i} \lambda q^{-n-1} w T \eta_n.
\end{equation*}

\bn
{\bf 7. Proofs of Theorems 1, 2 and 3}

\mn
Let us begin with some notation. For simplicity we write $x$ for 
the representation operator $\pi(x)$ of 
an algebra element $x\in\A$. Further, we shall omit the symbols 
$\lceil\D_0$ and $\lceil \D$ denoting the restrictions of the 
operators to $\D_0$ and $\D$, respectively. Moreover, we 
write simply $\Omega (.)$ instead of 
$\Omega_{\pi,F}$(.).

\mn
First let $\Gamma$ be the $3D$-calculus on $SU_q(2)$. We want to prove that 
the pair $(\pi,F)$ defined at the beginning of Section 5 is indeed a 
commutator representation of $\Gamma$. By Lemma 3, if suffices to show 
that $\Omega (x)=0$ for the six generators $x$ of the right ideal 
$\R_\Gamma$ listed by formula (\ref{3-36}). Computing the corresponding 
expressions of $\Omega (x)$ by using formula (\ref{l2-1}), we obtain the 
relations
\begin{gather}\label{5-61}
\Omega (q^2b^2)\equiv q^2 d^2 Fb^2+b^2Fd^2-(q^2+1) bd Fdb=0,\\
\label{5-62}
\Omega (c^2)\equiv q^2c^2Fa^2+a^2Fc^2-(q^2+1)ac Fca=0,\\
\label{5-63}
\Omega (qbc)\equiv-q^2 cd Fba+(q^2+1)bc Fbc-abFdc+qFbc+qbcF=0,\\
\label{5-64}
\Omega (q^2(a-1)b)\equiv q^2d^2 Fab-(q^2+1)bd Fbc+b^2Fcd-q^2dFb-qbdF+qbFd=0,\\
\label{5-65}
 \Omega ((a-1)c)\equiv -qcdFa^2+(q^2+1)bc Fca-ba Fc^2+Fac+qcFa-aFc=0,\\
\label{5-66}
\Omega ((q^2a+d-q^2-1))\equiv q^2dFa+aFd-qbFc-qcFb-(1+q^2)F=0.
\end{gather}
The relations (\ref{5-61}) and (\ref{5-64}) follows from (\ref{5-62}) and 
(\ref{5-65}), respectively, by applying the adjoint operation and 
using the equations (\ref{3-31}) and the fact that the operator $F$ is 
symmetric. Therefore it s sufficient to check 
(\ref{5-62}), (\ref{5-63}), (\ref{5-65}) and (\ref{5-66}). 
We omit these boring straightforward computations. In the course of 
these verifications the formulas (\ref{3-32})--(\ref{3-33}) and 
(\ref{4-44})--(\ref{4-45}) for the definition of the representation 
$\pi$ and of the operator $F$ and the relations 
(\ref{4-41})--(\ref{4-43}) are 
essentially used. Thus, $(\pi,F)$ is indeed a commutator representation 
of $\Gamma$. The admissibility of $(\pi,F)$ is obvious from its definition. 
This completes the proof of the first assertion of Theorem 1.

\bn
The next part of this section is devoted to the proofs of the second 
assertion of Theorem 1 and of Theorem 2. For this we suppose 
that $\Gamma$ is either the $3D$-calculus, the $4D_+$-calculus or 
the $4D_-$-calculus on $SU_q(2)$ and that $(\pi,F)$ is an arbitrary 
admissible commutator representation of $\Gamma$. Let $\D_0$, $\E$ and 
$\D$ be corresponding subspaces.

\bn
{\bf Lemma 5.} (i) {\it If $\Gamma$ is the $3D$-calculus, then we have} 
\begin{gather}\label{l5-1}
q^2c^2  Fa^2+a^2Fc^2-(q^2+1)acFca=0,\\
\label{l5-2}
qc^2Fa+qaFc^2-q^2c Fac - caFc=0,\\
\label{l5-3}
-q^2cdFba+(q^2+1) bcFbc -abFdc+ qFbc+qbcF=0,\\
\label{l5-4}
q^2dFa + aFd-qbFc - qcFb-(q^2+1)F=0.
\end{gather}
(ii) {\it If $\Gamma$ is the $4D_\pm$-calculus, then we have equations 
(\ref{l5-1}) and} 
\begin{gather}\label{l5-5}
a~Fc^2 -\epsilon ca Fc-q^2 c Fca +\epsilon qc^2 Fa=0,\\
\label{l5-9}
-q^2 dc Fa^2-ba Fc^2 -q^3c^2 Fab - qa^2Fcd + (q^2+1)bcFca
+q(q^2+1)acFbc +qFca + q^3caF=0.
\end{gather}

\bn
{\bf Proof.} (i): By (\ref{3-36}), the right ideal $\R_\Gamma$ associated 
with the $3D$-calculus contains the elements $c^2, bc ,q^2 a+d-(q^2+1)$ 
and $(a-1)c$. Therefore, by Lemma 2 we have $\Omega_{\pi,F}(c^2) =0 , 
\Omega_{\pi,F}(bc) =0$, $\Omega_{\pi,F} (q^2 a+d-(q^2+1))=0 $ and 
$\Omega_{\pi,F}((a-1)c)=0$. Computing these expressions by using 
(\ref{l2-1}) leads to the equations (\ref{l5-1}), (\ref{l5-3}), 
(\ref{l5-4}) and
\begin{equation}\label{l5-49}
(q^2+1)bcFca-qcdFa^2-baFc^2+Fac+qcFa-aFc = 0,
\end{equation}
respectively. It remains to derive equation (\ref{l5-2}). If we subtract 
equation $a$(\ref{l5-4})$c$ from $q^2$(\ref{l5-49}), we obtain
\begin{equation}\label{l5-7}
q^2bcFca + qacFbc - q^3cdFa^2-a^2Fdc + q^3cFa + aFc =0.
\end{equation}
Adding (\ref{l5-7}) and $qc$(\ref{l5-4})$a$  yields the equation 
\begin{equation}\label{l5-8}
(q^2+1)caFbc - q^2c^2Fba-a^2Fdc + qcaF+aFc-qcFa=0.
\end{equation}
Inserting the relation $a^2Fc^2=-q^2c^2Fa^2+(q^2+1)acFca$ by (\ref{l5-1})
into (\ref{l5-8})$c$ we finally get equation (\ref{l5-2}) as asserted.

\sn
(ii): Since the three elements $x=c^2, qc(a-d), z_\pm c$ belong to the 
right ideal associated with the $4D_\pm$-calculus (see (\ref{3-37})), 
the corresponding operators $\Omega_{\pi,F} (x)$ are zero by Lemma 2. 
This leads to the equations (\ref{l5-1}), (\ref{l5-9}) and 
\begin{align}\label{l5-10}
&-q^3 cdFa^2 + q^2 adFac + q^2 bcFca - q abFc^2 + q^2c^2 Fba
-qacFbc - qcaFda  \notag\\
&+a^2Fdc +\epsilon(q^4+1)
cFa -\epsilon (q^3+q^{-1})aFc=0, 
\end{align}
respectively. We still have to verify equation (\ref{l5-5}). 
Dividing (\ref{l5-9})--$q^2$(\ref{l5-10}) by $q^4+1$ and 
simplifying the terms by using the commutation rules of the matrix 
entries $a,b,c,d$, we obtain the equation 
\begin{equation}\label{l5-6}
(q^2+1) bcFca - baFc^2-q^2 dc Fa^2 + q Fca - q\epsilon a Fc+q^2 \epsilon 
c Fa=0.
\end{equation}
If we substitute $q^2c^2Fa^2 = (q^2+1)acFca-a^2 Fc^2$ 
(by (\ref{l5-1})) into $q^{-1}c$ (\ref{l5-6}), then 
equation (\ref{l5-5}) follows.\hfill\qed

\mn
Now we make use of the structure of the $\ast$-representation $\pi$ and 
of the admissibility of the pair $(\pi,F)$. We freely use the notation 
established above. Let $\D^n$ be the direct sum of domains 
$\D_k =\{ \eta_k: \eta \in \D_0 \}, k=0,\dots,n$, and let $\Hh^n$ 
denote the direct sum of subspaces $\Hh_k, k=0,\dots,n$, of  $\Hh$. Using 
essentially relation (\ref{l5-1}) and the 
fact that $\ker ~ a^k=\Hh_{k-1}$, a straightforward induction 
argument shows that the operator $F$ maps each space $\D^n$ into 
$\Hh^{n+1}$. This in turn implies that $F$ maps the subspace 
$\D = \rm{Lin} \{\D_n;n\in \bbbn_0\}$ into $\Hh = \oplus_n \Hh_n$. Since 
$F$ is symmetric, it follows that $F$ maps the domain $\E$ into $\G$. By 
assumption, $\E \oplus \D$ is a core for $F$. Therefore, the operator $F$ 
and hence all operators $\Omega_{\pi, F}(x), x \in \A$, leave 
the spaces $\G$ 
and $\Hh$ invariant. Using once more the facts that the operator $F$ is 
symmetric and that $F$ maps $\D^n$ into $\Hh^{n+1}$ it follows that the 
restriction of the operator $F$ to the dense linear subspace $\D$ of $\Hh$ 
is of the form 
\begin{equation}\label{l5-11}
F\eta_n= T_n \eta_{n-1} + R_n\eta_n+ 
T^\ast_{n+1} \eta_{n+1}, ~~\eta \in \D_0.
\end{equation}
Here $T_n$ and $R_n, n\in\bbbn_0$, are (possible unbounded) linear
operators on the Hilbert space $\Hh_0$ such that the domains of $T_n, R_n$ 
and $T^\ast_n$ contain $\D_0$ and $R_n$ is symmetric. Formula (\ref{l5-11}) 
will be essentially used in the sequel. For $n\in\bbbn$, we set 
\begin{equation}\label{l5-22}
E_n:=R_n-w R_{n-1} w^\ast.
\end{equation}
Inserting the formulas (\ref{3-33}) and (\ref{l5-11})) for the action of the 
operators $a,c$ and $F$ into (\ref{l5-1}) and comparing the expressions 
occuring in the $(n{-}2)$-th, $(n{-}1)$-th and $n$-th components, we obtain 
the recurrence relations
\begin{gather}\label{l5-12}
\lambda_{n+1}\lambda_nw^2 T_{n-1} + q^4 \lambda_n \lambda_{n-1} T_{n+1} 
w^2=(q+q^3)\lambda_{n+1} \lambda_{n-1} w T_n w, \\
\label{l5-13}
w^2 R_{n-1} + q^2 R_{n+1} w^2=(q^2+1) wR_n w,\\
\label{l5-14}
\lambda_{n+1}\lambda_n w^2 T^\ast_{n} + \lambda_{n+2} \lambda_{n+1} 
T^\ast_{n+2}  w^2=(q+q^{-1})\lambda^2_{n+1} w T^\ast_{n+1} w,
\end{gather}
respectively. Applying first the adjoint operation to (\ref{l5-14}), 
multiplying then by $w^2$ from the left and from the right, dividing by 
$\lambda_{n+1}$ and replacing finally $n$ by $n{-}1$, we get 
\begin{equation}\label{l5-15}
\lambda_{n-1} w^2 T_{n-1} + 
\lambda_{n+1} T_{n+1}w^2=(q+q^{-1})\lambda_n w T_n w.
\end{equation}
The equation $\lambda_{n-1}$(\ref{l5-12})--
$\lambda_n\lambda_{n+1}$(\ref{l5-15})
reads as 
$$
(q^4\lambda_n\lambda^2_{n-1} - \lambda_n \lambda^2_{n+1})T_{n+1} w^2= 
(q+q^3)\lambda_{n-1}^2 \lambda_{n+1} ((q+q^3)\lambda^2_{n-1} 
\lambda_{n+1}-(q+q^{-1}) \lambda^2_n \lambda_{n+1}) w T_n w.
$$
Since $\lambda^2_k= 1-q^{2k}$ by (\ref{3-35}), the latter yields 
$\lambda_n T_{n+1} w^2 = q^{-1} \lambda_{n+1} w T_nw$ and so 
\begin{equation}\label{l5-16}
q\lambda_n T_{n+1} = \lambda_{n+1} w T_n w^\ast.
\end{equation}
Note that the preceding formulas (\ref{l5-11}), (\ref{l5-13}) 
and
(\ref{l5-16}) are valid for both the $3D$-calculus and 
the $4D_\pm$-calculus, because they were derived only from 
equation (\ref{l5-1}) and this equation holds for all three calculi 
according to Lemma 5.

In order to continue the proof we first specify to the $3D$-calculus. 
Then, by Lemma 5(i), we have also equation (\ref{l5-2}). Inserting now 
(\ref{l5-11}) into (\ref{l5-2}) and comparing the $(n{-}1)$-th and $n$-th 
components, we get the relations
\begin{gather}\label{l5-17}
q^2 R_n w^2+w^2 R_{n-1}=q^2 w R_{n-1} w+ wR_nw,\\
\label{l5-18}
\lambda_{n+1} T^\ast_{n+1} w^2 + \lambda_n w^2 T^\ast_n = 
q\lambda_n w T^\ast_n w + q^{-1} \lambda_{n+1} w T^\ast_{n+1} w,
\end{gather}
respectively. Multiplying (\ref{l5-18}) by $w^\ast$ from the right, 
replacing $n$ by $n{-}1$, passing to the adjoint operators and 
finally applying formula (\ref{l5-16}), we derive 
\begin{equation}\label{l5-19}
\lambda_{n+1} T_n= \lambda_n T_{n+1}.
\end{equation}
Comparing (\ref{l5-16}) and (\ref{l5-19}) we conclude that 
$T_n = \lambda_n \lambda_1^{-1} T_1$ and $wT_nw^\ast = q T_n$. That is,
setting $T:= \lambda_1^{-1} T_1$, we have 
\begin{equation}\label{l5-20}
T_n = \lambda_n T ~~ {\rm and}~~ wTw^\ast  \eta = qT\eta,~ \eta \in \D_0.
\end{equation}

Next we investigate the diagonal terms $R_n$ of the operator $F$. 
First we note that in terms of the operators $E_n$ defined 
by (\ref{l5-22}) the equations (\ref{l5-13}) and (\ref{l5-17}) 
are reformulated as 
\begin{equation}\label{l5-21}
q^2 E_{n+1}=wE_n w^\ast ~~{\rm and}~~ q^2 E_n=w E_n w^\ast,
\end{equation}
respectively. In particular, we have $E_{n+1}=E_n$ for all $n$. 
If we compare the $n$-th components in (\ref{l5-3}), we get the relation 
\begin{equation}\label{l5-23}
q^{n+2}\lambda_n^2 wR_{n-1}w^\ast +(q^2+1) q^{2n+2} R_n + 
q^{n+2}\lambda^2_{n+1} w^\ast R_{n+1}w -2 q^{n+2} R_n =0.
\end{equation}
Putting the relations (\ref{l5-21}) into (\ref{l5-23}) we 
derive that $E_n=0$ for all $n$. Setting $R:=R_0$, the latter means that
\begin{equation}\label{l5-24}
R_n = w R_{n-1}w^\ast = w^n R w^{\ast n}.
\end{equation}
Further, comparing the $n$-th components in (\ref{l5-4}), we find that 
\begin{equation}\label{l5-25}
q^2 \lambda_n^2 R_{n-1} + \lambda_{n+1}^2 R_{n+1} 
+ q^{2n+2} w^\ast R_n w + q^{2n+2} w R_n w^\ast - (q^2+1)R_n = 0.
\end{equation}
Inserting the relation $E_n=0$ into (\ref{l5-25}), we obtain 
$R_{n+1} + q^2 R_{n-1}-(q^2+1)R_n =0.$ Because of (\ref{l5-24}), this means 
that
\begin{equation}\label{l5-26}
w^2Rw^{\ast 2} +q^2 R = (1+q^2)w R w^\ast.
\end{equation}

Finally, the restriction of the operator $F$ to the subspace $\E$ of the 
Hilbert space $\G$ is a symmetric linear operator, say $Q$. Since 
$b=c=0, a=v$ and $d=v^\ast$ on $\G$ by (\ref{3-32}), equation 
(\ref{l5-4}) reads as 
\begin{equation}\label{l5-27}
v^2Qv^{\ast 2} + q^2 Q = (1+q^2) vQv^\ast.
\end{equation} 

Summarizing the preceding, the formulas (\ref{l5-11}), (\ref{l5-20}), 
(\ref{l5-24}), (\ref{l5-26}) and (\ref{l5-27}) show that the operator 
$F$ has the required form. This completes the proof of the second assertion 
of Theorem 1. 

\mn
Now we turn to the $4D_\pm$-calculus and prove Theorem 2. To begin with, 
we compute the $(n{-}1)$-th and the $n$-th components of the expressions 
in equation (\ref{l5-5}). Comparing coefficients we derive 
\begin{gather}\label{l5-28}
R_nw^2 - q^{-1}\epsilon wR_n w-w R_{n-1} w + q^{-1} \epsilon w^2 R_{n-1}=0,\\
\label{l5-29}
\lambda_{n+1} T^\ast_{n+1} w^2 - \epsilon \lambda_{n+1} w T^\ast_{n+1} 
w - q\lambda_n w T^\ast_n w+q\epsilon \lambda_n w^2 T^\ast_n=0
\end{gather}
Applying the adjoint operation to (\ref{l5-29})$w^{\ast 2}$, we get 
\begin{equation}\label{l5-30}
\lambda_{n+1} T_{n+1} - \epsilon \lambda_{n+1}w T_{n+1} w^\ast
-q\lambda_n w T_n w^\ast + q\epsilon \lambda_n w^2 T_n w^{\ast 2}=0.
\end{equation}
Recall that formula (\ref{l5-16}) holds also for the $4D_\pm$-calculus,
because it was derived from formula (\ref{l5-1}). Inserting (\ref{l5-16})
into  (\ref{l5-30}), we derive that
\begin{equation}\label{l5-31}
T_{n+1}=\epsilon w T_{n+1} w^\ast.
\end{equation}
Combining the latter with (\ref{l5-16}), we get
\begin{equation}\label{l5-32}
q\lambda_n T_{n+1} =\epsilon \lambda_{n+1} T_n.
\end{equation}
Next we use equation (\ref{l5-9}) which holds by Lemma 5(ii).
Computing the $(n{-}1)$-th components of 
(\ref{l5-9}), we obtain the relation
\begin{align}\label{l5-34}
&-q^n \lambda_n \lambda_{n-1}^2 w R_{n-2} + 
q^{3n} \lambda_n w^\ast R_n w^2 +q^{3n+2} \lambda_n w^2 R_{n-1} w^\ast - 
q^{n+2} \lambda_{n+1}^2 \lambda_n R_{n+1} w\notag\\ &-(q^2+1)q^{3n-2} 
\lambda_n R_{n-1} w
-(q^2+1)q^{3n+2} \lambda_n w R_n +q^n \lambda_n R_{n-1}w 
+q^{n+2} \lambda_n w R_n = 0.
\end{align}
In terms of the 
operator $E_n$, the formulas (\ref{l5-13}) and (\ref{l5-28}) can
be written as 
\begin{equation}\label{l5-35}
q^2 E_{n+1} =w E_nw^\ast ~{\rm and}~ E_n=\epsilon q^{-1} w E_n w^\ast
\end{equation}
respectively. Inserting these formulas into (\ref{l5-34})$w^\ast$, a 
lengthy computation shows that $E_n=0$. That is, setting $R:=R_0$, we have 
\begin{equation}\label{81}
R_n =w R_{n-1} w^\ast = w^n R w^{\ast 2}\mbox{ for } n\in\bbbn_0.
\end {equation}
Using  the formulas (\ref{l5-10}), (\ref{l5-6}) and (\ref{81}) established above, 
we compute
\begin{align}\label{82}
\Omega (a)\eta_n &= (\epsilon q-1) T_n \eta_{n-1}+ (R_{n-1}-R_n)\eta_n + 
(\epsilon q-1) T^\ast_{n+1} \eta_{n+1},\\
\label{83}
\Omega (b) \eta_n &= - \lambda q^n \lambda^{-1}_n T_n w^\ast 
\eta_n \equiv \lambda^{-1}_q \lambda^{-1}_n T_n b \eta_n,\\
\label{84}
\Omega (c) \eta_n &=-\lambda q^n \lambda^{-1}_n w T^\ast_n \eta_n 
\equiv -\lambda\lambda^{-1}_n c T^\ast_n \eta_n,\\
\label{85}
\Omega (d) \eta_n&= (\epsilon q^{-1} -1)T_u \eta_{n-1} + 
(R_{n+1}- R_n)\eta_n+(\epsilon q^{-1} -1)T^\ast_{n+1} \eta_{n+1},
\end{align}
for $\eta \in\D_0$ and $n\in \bbbn$. In particular, we get 
\begin{equation}\label{86}
\Omega (a+\epsilon q d)\eta_n= (\epsilon q R_{n+1}-(\epsilon q+1) R_n+R_{n-1})\eta_n.
\end{equation}
Put $\Omega_j=\rho (-{\rm i}\omega_j)$. Since $\rho$ is bimodule homomorphism, the commutation relations between the 1-forms $\omega_j$ and the generators 
$a,b,c,d$ remain valid if $\omega_j$ is replaced by $\Omega_j$. In 
particular, the relation $\omega_2 a=\epsilon a\omega_2 + 
\epsilon \lambda^2 q^{-1}\omega_4$ yields
\begin{equation}\label{87}
\Omega_2 a=\epsilon a\Omega_2+\epsilon\lambda^2 q^{-1}b\Omega_4.
\end{equation}
Since $\Omega (a+\epsilon q d)=(1-q^2)(\epsilon q^{-3}-1)\Omega_4$ by 
(\ref{78}) 
and (\ref{79}), it follows from (\ref{86}) that 
$\Omega_4$: $\Hh_n\rightarrow \Hh_n$ and so 
$b\Omega_4:\Hh_n\rightarrow \Hh_n$. Further, by (\ref{83}) the 
operators $\Omega_2a$ and $\epsilon a\Omega_2$ map 
$\Hh_n$ into $\Hh_{n-1}$. Because of (\ref{87}) this implies that 
$\Omega_4=0$ on each space $\Hh_n$ and so on $\Hh$.
From the relation $\omega_2c=\epsilon c\omega_2+\epsilon\lambda^2 q^{-1} 
d\omega_4$ we 
get $\Omega_2c=\epsilon c\Omega_2+\epsilon \lambda^2 q^{-1} d\Omega_4.$ 
Since $c=0$ and 
$d$ is unitary on $\G$, this implies that $\Omega_4=0$ on $\G$.
 
Thus we have shown that $\Omega_4=0$ on $\G\oplus \Hh$ which yields that 
$\rho (\omega_{\Gamma_\pm}(a+\epsilon qd))=0$. Therefore, by Lemma 3 the 
FODC homomorphism 
$\rho$ passes to the quotient FODC $\Gamma_{\pm,3}$. This completes the 
proof of Theorem 2.

\bn
{\bf Remark 3.} The assertions of Theorems 1 and  2 could be also 
derived from the commutation relations and the involution properties 
of the calculi thus avoiding the use of the right ideals. 
We prefered to give the above proof 
because it emphasizes the role of the corresponding right ideals and it 
needs only very few generators of the right ideals. 

\mn
{\bf Remark 4.} By adding a few lines to the preceding 
arguments one gets a complete description of all admissible 
commutator representations of the quotient $\ast$-FODC 
$\Gamma_ {\pm,3}$. In this remark we briefly derive this 
result. Since $\Omega_4=0$ and hence $\Omega (a+\epsilon qd)=0$, it 
follows from (\ref{81}), (\ref{82}) and (\ref{85}) that 
\begin{equation}\label{90}
\omega^2R\omega^{\ast 2}+\epsilon q^{-1} R=(1+\epsilon q^{-1}) 
w Rw^\ast\mbox{ on } \D_0.
\end{equation}
Set $T:=\lambda^{-1}_1 T_1$. From (\ref{l5-31}) and (\ref{l5-32}) we get 
$T_n=(\epsilon q)^{1-n} \lambda_nT$ and 
\begin{equation}\label{91}
wTw^\ast=\epsilon T\mbox{ on } \D_0.
\end{equation}
Let $Q$ denote the restriction of $F$ to the domain $\E$ in the 
subspace $\G$. Since $b=0$ on $\G$, we have 
$\Omega (b)\equiv dFb-q^{-1} bFd=0$ and so $\Omega_2=0$ on $\G$. 
Hence the relation $\Omega_1 d=\epsilon q^{-1} d\Omega_1+\epsilon c\Omega_2$ 
implies that $\Omega_1 v^\ast=\epsilon q^{-1} v^\ast\Omega_1$ on $\G$. 
Because $\Omega_4=0$ as shown in the 
above proof, we have $\Omega (a)= (\epsilon q-1)\Omega_1$ by 
(\ref{78}). On $\G$ we have $\Omega (a)\equiv dFa-q^{-1} bFc-F=v^\ast Qv-Q$. 
Inserting the latter into the relation 
$\Omega (a)v^\ast=\epsilon q^{-1} v^\ast \Omega (a)$, we finally obtain that
\begin{equation}\label{92}
v^2Qv^{\ast 2}+\epsilon q^{-1} Q=(1+\epsilon q^{-1})vQv^\ast\mbox{ on } \E.
\end{equation}
The symmetric operator $F$ now acts as
\begin{equation}\label{93}
F\eta_n=(\epsilon q)^{1-n} \lambda_n T \eta_{n-1} + 
w^n Rw^{\ast n} \eta_n+ (\epsilon q)^{-n} \lambda_{n+1} 
T^\ast \eta_{n+1}, \eta\in \D_0, \mbox{ and  } F\eta=Q\eta,\eta\in\E.
\end{equation}

Conversely, let $T$ be a linear operator and $R$ a symmetric linear operator 
on a dense domain $\D_0$ of $\Hh$ and let $Q$ be a symmetric linear 
operator on a dense domain $\E$ of $\G$  such that 
$\D_0\subseteq\D (T^\ast), w\D_0=\D_0$ and $v\E=\E$. If the relations 
(\ref{90}), (\ref{91}) and (\ref{92}) are fulfilled, then the pair 
$(\pi,F)$ with $F$ defined by (\ref{93}) is an admissible commutator 
representation of the $\ast$-FODC $\Gamma_{\pm,3}$.

\bn

Now let us turn to the proof of Theorem 3. The key of the proof is 
the following simple lemma.

\mn
{\bf Lemma 6.} {\it Let $u$ be a unitary operator and let $A$ be a bounded 
linear operator on a Hilbert space. If $uAu^\ast=\alpha A$ for some 
$\alpha\in\bbbr, |\alpha |\ne 1,$ then $A=0$.}

\bn
{\bf Proof.} Then we also have $uA^\ast u^\ast = \alpha A^\ast$. So we may 
assume that $A$ is self-adjoint. The relation $uAu^\ast=\alpha A$ implies 
that the spectrum of $A$ is invariant under multiplication by $\alpha$ and 
$\alpha^{-1}$. Since $|\alpha |\ne 1$ and $A$ is bounded, this is only 
possible if $A=0$.\hfill \qed

\mn
We carry out the proof for the $4D_\pm$-calculus. The case of the 
$3D$-calculus is much simpler and follows easily from the relations 
(\ref{70}) and (\ref{71}) and Lemma 6.

By assumption all operators $\dif_{\pi,F}(x), x\in\A$, are bounded. 
Hence the four operators $\Omega_1, \Omega_2, \Omega_3, \Omega_4$ 
are also bounded. 
Recall that $\Omega^\ast_1=\Omega_1, \Omega^\ast_2=\Omega_3$ and 
$\Omega^\ast_4=\Omega_4$ by (\ref{89}). From these facts and the 
commutation relations of the $4D_\pm$-calculus it is clear that the 
operators $\Omega_j$ leave the spaces $\G=\kere~b=\kere~c$ and $\Hh$ 
invariant.

We shall show that all $\Omega_j=0$ on $\Hh$ for $j=1,{\dots}, 4$. First 
note that $\Omega_4:\Hh_n\rightarrow\Hh_n$, because of the relation 
$\Omega_4 bc=bc\Omega_4$ and the fact that $\Hh_n=\ker ~(bc+q^{2n+1} I).$
Since $\Omega_4 b=\epsilon q b\Omega_4$, we therefore obtain that 
$w(\Omega_4\lceil \Hh_n)w^\ast = \epsilon q (\Omega_4\lceil\Hh_n)$, so 
that $\Omega_4\lceil \Hh_n=0$ by Lemma 6 and hence 
$\Omega_4=0$ on $\Hh$. Consequently we have $\Omega_3 bc =bc\Omega_3$ 
which implies that 
$\Omega_3:\Hh_n\rightarrow\Hh_n$. Since $\Omega_4=0$, we have the 
two relations
\begin{gather}\label{94}
\Omega_1 a=\epsilon q a\Omega_1+\epsilon b\Omega_3,\\
\label{95}
\Omega_1 b=\epsilon q^{-1} b\Omega_1+\epsilon a\Omega_2,
\end{gather}
Using the fact that $\Omega_3:\Hh^n\rightarrow\Hh^n$ it follows from 
(\ref{95}) by induction that $\Omega_1:\Hh^n\rightarrow \Hh^{n+1}$. 
Because $\Omega^\ast_1=-\Omega_1,$ we have 
$\Omega_1:\Hh_n\rightarrow\Hh_{n-1}\oplus \Hh_n \oplus \Hh_{n+1}$. 
Hence $\Omega_1$ 
is of the form
$$
\Omega_1\eta_n= A_n\eta_{n-1}+B_n\eta_n+ A^\ast_{n+1} 
\eta_{n+1},
\eta\in\Hh_0,
$$
where $A_n$ and $B_n$ are bounded linear operators on $\Hh_0$. Inserting 
this expression into (\ref{95}) and comparing the $n$-th components, we 
get $qB_n w^\ast=\epsilon w^\ast B$. Thus, $wB_nw^\ast=\epsilon q^{-1}B_n$ 
and hence $B_n=0$ by Lemma 6. Comparing the $(n{-}2)$-th components 
in (\ref{94}), we obtain the equation 
$\lambda_n A_{n-1}=\epsilon q\lambda_{n-1} A_n$, so that we have
\begin{equation}\label{98}
A_n=(\epsilon q)^{1-n} \lambda_n\lambda^{-1}_1 A_1
\end{equation}
Since $\parallel A_n\parallel~ \le~ \parallel\Omega_1\parallel$ and 
$q^{-n}\lambda_n\rightarrow +\infty$ if $n\rightarrow \infty$ 
(recall that $0<q<1$), we conclude from (\ref{98}) 
that $A_n=0$ for all $n\in\bbbn$. 
Thus  $\Omega_1=0$ on $\Hh$. Applying once more equations (\ref{94}) 
and (\ref{95}), we see that $\Omega_2=\Omega_3=0$ on $\Hh$.

A much simpler reasoning shows that the operators $\Omega_j$ are also 
zero on $\G$. Since the four $1$-forms $\omega_j, j=1,{\dots},4$, 
generate 
the $4D_\pm$-calculus as a left $\A$-module, it follows that 
$\dif_{\pi,F}(x)=0$ for all $x\in \A$. This completes the 
proof of Theorem 3. \hfill\qed

\bn
{\bf Acknowledgement:} This work was completed during the authors visit 
at the University of Fukuoka in September 1998. I am very grateful to  
Professors A. Inoue and H. Kurose and Dr. H. Ogi for their warm hospitality 
and for the excellent working atmosphere.

\sn


\end{document}